\documentclass[final,3p]{elsarticle}

\usepackage{stmaryrd} %for DG jump brackets	
\usepackage{geometry,graphicx}
\usepackage[normalem]{ulem}
\usepackage{color}
\usepackage[colorlinks=true,bookmarks,bookmarksopen,bookmarksdepth=2]{hyperref}
\usepackage[all]{hypcap}
\usepackage[T1]{fontenc} 							
\usepackage{amssymb,amsmath}%,amsbsy}%amsthm		
\usepackage{empheq}							
\usepackage{bm} % for \boldsymbol{}
\usepackage{thmbox}
\usepackage{shadethm}		
\usepackage{scalerel} %for \fakebold
\usepackage{nicefrac}								
\usepackage{siunitx}  		
\usepackage{booktabs} %for nice rules in tables
\usepackage{subfigure}		
\usepackage{mathtools}	% for e.g. \coloneqq
\usepackage{lipsum}
\usepackage{pgfplots}

\usepackage{titlesec}				% for formatting of CAPTIONS	 

\titleformat{\section}{\fontsize{11}{17}\bfseries}{\thesection}{1em}{}
\titleformat{\subsection}{\fontsize{10}{17}\bfseries\itshape}{\thesubsection}{1em}{}		

\newcommand\redsout{\bgroup\markoverwith{\textcolor{red}{\rule[0.5ex]{2pt}{0.5pt}}}\ULon}

\newcommand\bluesout{\bgroup\markoverwith{\textcolor{blue}{\rule[0.5ex]{2pt}{0.5pt}}}\ULon}

\newlength\bshft
	\def\fakebold#1{\ThisStyle{\ooalign{$\SavedStyle#1$\cr%
  	\kern-\bshft$\SavedStyle#1$\cr%
  	\kern\bshft$\SavedStyle#1$}}}

\newcommand{\R}{\mathbb{R}}
\newcommand{\N}{\mathbb{N}}

\newcommand{\dvg}{{ \mathrm{div} }}

\newcommand{\im}{\mathrm{i}}
\newcommand{\intrm}{{ \mathrm{int} }}
\newcommand{\rel}{{ \mathrm{rel} }}
\newcommand{\rms}{{ \mathrm{RMS} }}
\newcommand{\OMEGA}{{ \rb{\Omega} }}
\newcommand{\tbar}{{ \bar{t} }}

\newcommand\Rey{\mbox{\textit{Re}}}  

\DeclareMathOperator{\ip}{{\boldsymbol{\cdot}}}

\DeclareMathOperator{\DIV}{\nabla\,\boldsymbol{\cdot}}

\DeclareMathOperator{\crit}{{ crit }}

\newcommand{\uu}{{ \boldsymbol{u} }}
\newcommand{\vv}{{ \boldsymbol{v} }}

\newcommand{\xx}{{ \boldsymbol{x} }}

\newcommand{\nn}{{ \boldsymbol{n} }}
\newcommand{\kk}{{ \boldsymbol{k} }}

\newcommand{\zero}{{ \boldsymbol{0} }}
\newcommand{\tend}{{ T }}

\newcommand{\drm}{{ \mathrm{d} }}
\newcommand{\dxx}{{ \,\drm\xx }}
\newcommand{\dx}{{ \,\drm x }}
\newcommand{\dy}{{ \,\drm y }}
\newcommand{\ds}{{\,\drm\boldsymbol{s}}}
\newcommand{\dtau}{{ \,\drm\tau }}

\newcommand{\eps}{{ \varepsilon }}

\newcommand{\lavg}{{ \big\{\hspace{-0.99ex}\big\{ }}						
\newcommand{\ravg}{{ \big\}\hspace{-0.99ex}\big\} }}		
\newcommand{\ljmp}{ \left\llbracket }	% no double {{ }}									
\newcommand{\rjmp}{ \right\rrbracket }	% no double {{ }}								
\newcommand\jmp[1]{{ \ljmp#1\rjmp }}										
\newcommand\avg[1]{{ \lavg#1\ravg }}

\newcommand{\Kin}{{ \mathcal{K} }}
\newcommand{\Ens}{{ \mathcal{E} }}
\newcommand{\Pal}{{ \mathcal{P} }}
\newcommand{\Quo}{{ \mathcal{Q} }}

\newcommand\Lp[2]{{ L^{#1}{#2} }} 
\newcommand\LP[2]{{ \boldsymbol{L}^{#1}{#2} }}

\newcommand\HM[2]{{ \boldsymbol{H}^{#1}{#2} }}

\newcommand\HDIV{{ \boldsymbol{H}{\rb{\dvg}} }}
\newcommand{\Hdiv}{{ \boldsymbol{H}{\rb{\dvg;\Omega}} }}

\newcommand{\T}{{ \mathcal{T}_h }}

\newcommand{\Fi}{{ \mathcal{F}_h^i }}

\newcommand\rb[1]{{ \left(#1\right) }}
\newcommand\sqb[1]{{ \left[ #1 \right] }}
\newcommand\rsb[1]{{ \left(#1\right] }}
\newcommand\set[1]{{ \left\{ #1 \right\} }}
\newcommand\bra[1]{{ \langle #1 \rangle }}
\newcommand\abs[1]{{ \left\lvert#1\right\rvert }}
\newcommand\norm[1]{ \left\lVert#1\right\rVert }

\newcommand\nf[2]{{ \nicefrac{#1}{#2} }}

\newcommand{\tripnorm}[1]{{\left\vert\kern-\nulldelimiterspace\left\vert\kern-\nulldelimiterspace\left\vert #1
	\right\vert\kern-\nulldelimiterspace\right\vert\kern-\nulldelimiterspace\right\vert}}
	
\newcommand{\otoprule}{\midrule[\heavyrulewidth]}

\newcommand\restr[2]{{												
	\left.\kern-\nulldelimiterspace									
	#1
	\vphantom{\big|}
	\right|_{#2}
	}}
	
\newcommand{\goodgap}{%
	\hspace{0.01\subfigtopskip}
	\hspace{0.01\subfigbottomskip}
	}

\newtheorem[style=S,underline=true,bodystyle=\normalsize\noindent]{thmDef}{\textsc{Definition}}[section]
\newtheorem[style=S,cut=false]{thmCor}[thmDef]{\textsc{Corollary}}
\newtheorem[style=S,cut=false,headstyle=\normalsize\bfseries\boldmath####1~####2]{thmLem}[thmDef]{\textsc{Lemma}}

\newshadetheorem{thm}[thmDef]{\textsc{Theorem}}

\newenvironment{thmProof}
                [1][\unskip]
                { \begin{example}[\normalsize \textsc{Proof #1}] \normalsize}
                { $\hfill\blacksquare$ \end{example} }	                
                
\newenvironment{thmRem}
                [0]
                { \refstepcounter{thmDef} \begin{example}[\normalsize\textsc{Remark} \thesection.\arabic{thmDef}]  \normalsize}
                { \end{example} }	
                
\newtheorem[style=S,underline=true,bodystyle=\noindent,cut=false]{thmAss}{\small\textsc{Assumption}}

\journal{`Computers \& Mathematics with Applications' (accepted: October 20, 2018)}
\definecolor{mediumblue}{RGB}{0,0,205}
\definecolor{forestgreen}{RGB}{34,139,34}
\definecolor{darkred}{RGB}{200,0,0}

\begin{document}

\hypersetup{
  linkcolor=darkred,
  urlcolor=forestgreen,
  citecolor=mediumblue
}

%------------------------------------------------------------------------------------------------
\begin{frontmatter}


\title{On reference solutions and the sensitivity of the \\ 2D Kelvin--Helmholtz instability problem} 

\author[goe]{Philipp W.\ Schroeder\fnref{fn1}}
\ead{p.schroeder@math.uni-goettingen.de}
\address[goe]{Institute for Numerical and Applied Mathematics, Georg-August-Universit{\"a}t G{\"o}ttingen, 37083 G\"ottingen, Germany}
\fntext[fn1]{ORCID: \url{https://orcid.org/0000-0001-7644-4693}}

\author[wias,fub]{Volker John\fnref{fn2}}
\ead{john@wias-berlin.de}
\address[wias]{Weierstrass Institute for Applied Analysis and Stochastics, 10117 Berlin, Germany}
\address[fub]{Department of Mathematics and Computer Science, Freie Universit{\"a}t Berlin, 14195 Berlin, Germany}
\fntext[fn2]{ORCID: \url{https://orcid.org/0000-0002-2711-4409}}  

\author[tuw]{Philip L.\ Lederer} 
\ead{philip.lederer@tuwien.ac.at}
\address[tuw]{Institute for Analysis and Scientific Computing, TU Wien, 1040 Wien, Austria} 

\author[goe]{Christoph Lehrenfeld\corref{cor1}\fnref{fn4}}
\ead{lehrenfeld@math.uni-goettingen.de}
\cortext[cor1]{Corresponding author}
\fntext[fn4]{ORCID: \url{https://orcid.org/0000-0003-0170-8468}}

\author[goe]{Gert Lube}
\ead{lube@math.uni-goettingen.de}

\author[tuw]{Joachim Sch{\"o}berl} 
\ead{joachim.schoeberl@tuwien.ac.at}

%------------------------------------------------------------------------------------------------
\begin{abstract}

Two-dimensional Kelvin--Helmholtz instability problems are popular examples for assessing 
discretizations for incompressible flows at high Reynolds number.
Unfortunately, the results in the literature differ considerably. 
This paper presents computational studies of a Kelvin--Helmholtz instability problem with high order divergence-free finite element methods. 
Reference results in several quantities of interest are obtained for three different Reynolds numbers up to the beginning of the final vortex pairing. 
A mesh-independent prediction of the final pairing is not achieved due to the sensitivity of the considered problem with respect to small perturbations. 
A theoretical explanation of this sensitivity to small perturbations is provided based on the theory of self-organization of 2D turbulence. 
Possible sources of perturbations that arise in almost any numerical simulation are discussed. 

\end{abstract}
%------------------------------------------------------------------------------------------------

%------------------------------------------------------------------------------------------------
%-------------------------------------------------------------------------------------------------

\hspace{-5.25mm}Publisher's version: DOI \url{https://doi.org/10.1016/j.camwa.2018.10.030}  \\
\copyright~2018. This manuscript version is made available under the CC BY-NC-ND 4.0 license: \\
\url{https://creativecommons.org/licenses/by-nc-nd/4.0/}\\

\end{frontmatter}

%------------------------------------------------------------------------------------------------
%----------------------------------------INTRODUCTION--------------------------------------------
%------------------------------------------------------------------------------------------------
\section{Introduction}
\label{sec:Introduction}

Before starting to read, we recommend at first to search the internet for `Kelvin Helmholtz instability nature'.
Mostly, this search will reveal photographs which show the physical phenomenon in the context of cloud formation.
Take a closer look at them and, hopefully, appreciate the beauty of it!
Note, however, that in nature the occurrence of any kind of perturbation is frequent and absolutely normal.
This fact will be a recurring topic throughout this paper.

Good benchmark problems are necessary for assessing numerical schemes. 
For high Reynolds number incompressible flows in three dimensions, i.e., turbulent flows, there are a number of commonly used problems, like the isotropic turbulence problem, where the behavior of the energy spectrum is known from the K41 theory \cite{Kol41,Kol41a}, turbulent channel flows at several Reynolds numbers \cite{MoserKimMansour99}, or the turbulent flow around a cylinder \cite{RodiEtAl97}. 
It is well known that real turbulence possesses some features that are inherently three-dimensional, like vortex stretching. 
But in our opinion, also good two-dimensional benchmark problems for high Reynolds number flows are of interest.
There is the same principal difficulty in 2D and 3D: the flow possesses important scales that, depending on the particular problem, can be difficult to resolve. 
A stabilization technique (turbulence  model in 3D) has to take into account the impact of small and unresolved scales on the computed solution. 
A frequently used approach for assessing proposals for discretizations of high Reynolds number flows consists in studying first 2D examples, because of the easier implementation and usually much shorter computing times compared with 3D, before applying them in turbulent 3D simulations.
However, note that it is an open question to which extent the measures and metrics with which we assess numerical results in this work can be translated to 3D.
A crucial property of the Kelvin--Helmholtz instability problem considered in this work turns out to be its high sensitivity which allows one to critically assess the discretization accuracy of discrete flow solvers.

In the literature, the most often used example for a high Reynolds number flow in two dimensions seems to be the Kelvin--Helmholtz instability or mixing layer problem defined in \cite{LesieurEtAl88}.
Starting from a noisy initial condition, small vortices arise which then pair to larger and larger vortices until finally one vortex remains. 
Note that such a behavior of energy transfer from small to large scales is characteristic for (high Reynolds number) two-dimensional flows \textemdash{} three-dimensional flows usually cannot reorganize themselves into large structures. 
The Kelvin--Helmholtz instability example presents a richness of flow scales and an interesting temporal evolution of the flow field. 
Furthermore, the mixing layer problem possesses the classical features of deterministic chaos inherent to the Navier--Stokes problem.
For these reasons, we think that it fits very well for assessing numerical schemes for 2D turbulent incompressible flow simulations. 

A review of numerical studies for the Kelvin--Helmholtz instability problem will be provided in Section~\ref{sec:QuantitiesOfInterest}. 
Checking the presented results, one can see that all of them are qualitatively correct in the sense that vortex pairings up to one big vortex are predicted. 
However, quantitatively, the results are often considerably different, in particular with respect to the times at which the pairings occur. 
Maybe most notably, one can find results where the pairing from two vortices to one vortex follows very shortly after the previous pairing, e.g., in \cite{John05,YangBadiaCodina16}. 
In other numerical studies, there is a comparatively long interval between those pairings, e.g., in \cite{GravemeierEtAl05,SchroederLube18}. 
Since there are no reference results available, it is not clear which behavior is the correct one.
This last vortex pairing is crucial in this work and it turns out that it cannot be predicted reliably due to the sensitivity of the problem.

The derivation of reference results for the 2D Kelvin--Helmholtz instability problem would be very helpful for defining a good benchmark problem for high Reynolds number 2D flows. 
To compute such results, modern numerical methods and the nowadays available computational power should be used. 
In the numerical studies presented in this paper, direct numerical simulations (DNS) with  $\HDIV$-conforming discrete velocity spaces of polynomial degree $8$ were performed. 
Being $\HDIV$-conforming means that the discrete velocity field is divergence-free in the sense of $\Lp{2}{\OMEGA}$, a property which is not given by most standard discretizations. 
The negative impacts of not being $\HDIV$-conforming on the results computed with finite element methods have been recognized only recently, see the review paper \cite{JLMNR17}. 

As a first contribution of this paper, reference results for the 2D Kelvin--Helmholtz instability problem for several Reynolds numbers up to the evolution of two vortices at around 200 (scaled) time units are provided. 
The second result of interest is that even with the used state-of-the-art methods a conclusive prediction of the final pairing to one vortex is not yet possible.
To be concrete, even among the simulations with the highest resolutions, the time interval for the final pairing is still somewhat different.
This is not only due to a very high demand on resolution in space and time, but foremost due to the high sensitivity of the problem to small perturbations that are almost unavoidable in numerical simulations.
Thirdly, we apply the theory for self-organization of 2D turbulent flows to give a consistent explanation, on the continuous level, for this high sensitivity.
The computational results are made available for the community at \url{https://ngsolve.org/kh-benchmark} \cite{KHpage}.

\emph{Organization of the article:}
In Section~\ref{sec:KHInstabilityProblem}, the general setting for the Kelvin--Helmholtz instability problem which is considered in this paper is explained.
Also, frequently evaluated quantities of interest are introduced and discussed.
Section~\ref{sec:SelfOrga2DTurbulence} deals with the aspect of self-organization in turbulent two-dimensional incompressible flows. 
Based on a theory by Van Groesen \cite{VanGroesen88}, it is intended to raise the readers' awareness to the fact that already on the continuous level, such flow systems are very sensitive with respect to perturbations.
In Section~\ref{sec:DivFreeHDG}, we describe the numerical method that was used for the simulations. 
Both time and space discretization are addressed.
Then, Section~\ref{sec:ComputationalStudies} presents numerical results for three different Reynolds numbers.
Several quantities of interest are evaluated and discussed there.
In order to even more emphasize the difficulties in computing reference solutions for this problem, Section~\ref{sec:Perturbations} shows how small perturbations stemming from various sources in the process of discretization can significantly change the outcome of the simulations.
We summarize and conclude in Section~\ref{sec:Conclusions}.
Lastly, with the intention of making it easy for the community to compare their own results with ours, the Appendix~\ref{sec:Appendix} explains how  the data of our reference results can be accessed. 

%------------------------------------------------------------------------------------------------
%----------------------------------KH INSTABILITY PROBLEM----------------------------------------
%------------------------------------------------------------------------------------------------
\section{The Kelvin--Helmholtz instability problem}
\label{sec:KHInstabilityProblem}

The most frequently found setting for a Kelvin--Helmholtz instability problem is the evolution/dissipation of an initial condition in a viscous incompressible Navier--Stokes flow. 
Due to the fact that no body forces are present, the whole motion is thus driven only by the initial condition. 
For the kinematic viscosity $\nu >0$, the time-dependent incompressible Navier--Stokes problem with vanishing source term reads:
\begin{subequations}\label{eq:TINS}
\begin{empheq}[left=\empheqlbrace]{alignat=2}  
  	\partial_t \uu - \nu \Delta \uu + \rb{\uu \ip \nabla} \uu + \nabla p &= \zero, 
  		\qquad &&\text{in }\rsb{0,\tend} \times \Omega, \\
  	\DIV \uu &= 0, 
  		\qquad &&\text{in }\rsb{0,\tend} \times \Omega, \\
  	\uu\rb{0,\cdot} &= \uu_0\rb{\cdot }, 
  		\qquad &&\text{in } \Omega.	
\end{empheq}	
\end{subequations}

%------------------------------------PROBLEM STATEMENT-------------------------------------------
\subsection{Problem statement}
\label{sec:ProblemStatement}

The computational domain for the Kelvin--Helmholtz instability problem is a square. 
In this paper, $\Omega = \rb{0,1}^2$ is considered, but one finds in the literature also other setups, e.g., with $\Omega = \rb{0,2\pi}^2$ as in  \cite{SchneiderFarge00} or $\Omega = \rb{-1,1}^2$ as in \cite{John05}.
At $x=0$ and $x=1$, periodic boundary conditions are used, mimicking in this way an infinite extension in horizontal direction. 
At $y=0$ and $y=1$, free-slip boundary conditions are prescribed. 
The initial condition is given by 
\begin{align}\label{eq:ini_cond}
	\uu_0\rb{x,y} 
		= \begin{bmatrix} 
			u_\infty \tanh \rb{ \frac{2y-1}{\delta_0} } \\
			0 
		\end{bmatrix}
		+ c_n
		\begin{bmatrix} 
			\partial_y \psi\rb{x,y} \\ 
			- \partial_x \psi\rb{x,y} 
		\end{bmatrix}
\end{align}
with corresponding stream function
\begin{align*}
	\psi\rb{x,y} 
		= u_\infty\exp\rb{-\frac{\rb{y-0.5}^2}{\delta_0^2}}
			\sqb{ \cos\rb{8\pi x} + \cos\rb{20\pi x} }.
\end{align*}
Here, $\delta_0=1/28$ denotes the initial vorticity thickness, $u_\infty=1$ is a reference velocity and $c_n=\num{E-3}$ is a scaling/noise factor.
Note that frequently the $\cos\rb{20 \pi x}$ term is not included in the literature.

The principal behavior of the flow is as follows, compare \cite[Sec.\ 3.4.1]{Lesieur08} or, originally \cite{Michalke64}. 
The perturbations prescribed in the right-hand side term of the initial condition \eqref{eq:ini_cond} are amplified such that vortices develop. 
Here, the most amplified mode corresponds to the longitudinal wavelength $\lambda_a = 7 \delta_0$. 
In particular, $n\in \N$ primary vortices develop in a domain with length $n\lambda_a$ in horizontal direction. 
Hence, due to choosing $\delta_0 = 1/28$, in the numerical simulations, $n=4$, i.e.,  so-called `4-eddy calculations' \cite{LesieurEtAl88}, were performed.

The Reynolds number $\Rey$ of the Kelvin--Helmholtz instability flow is usually calculated on the bases of the characteristic length scale $\delta_0$ and the characteristic velocity scale $u_\infty$, i.e., $\Rey=  \delta_0 u_\infty/\nu = 1/\rb{28\nu}$. 
The numerical simulations will study different values of $\nu$ such that we consider $\Rey\in\set{\num{100},\num{1000},\num{10000}}$. 
For the simulations and their evaluation, the time unit $\tbar = \delta_0/u_\infty$ is introduced.

%---------------------------------QUANTITIES OF INTEREST-----------------------------------------
\subsection{Quantities of interest}
\label{sec:QuantitiesOfInterest}

For the definition of benchmark problems, it is necessary to define appropriate quantities of interest. On the one hand, these quantities should be of some physical importance. 
But on the other hand, it is of advantage for the benchmark problem to be accepted by the community if the implementation of their calculation can be done with affordable effort. 
This section reviews the quantities of interest that were considered in the literature and points out those quantities of interest that are studied in this paper. 

%---------
\emph{Mean streamwise velocity profiles with RMS profiles.}
Given a point in time $t^*$, the mean streamwise velocity profile is given by 
\begin{align*}
	\frac{\int_{0}^1 u_1\rb{t^*, x, y} \dx}{\int_{0}^1 \dx} 
		=  \int_{0}^1 u_1\rb{t^*, x, y} \dx, 
		\quad y\in \sqb{0,1}.	
\end{align*}
In practice, in particular if equidistant meshes are used, the integral definition is replaced by an arithmetic average
\begin{align}\label{eq:mean_velo}
	\bra{u_1}\rb{t^*, y} 
		= \frac{1}{N_x} \sum_{i=1}^{N_x} u_1\rb{t^*, x_i, y}, 
		\quad y\in \sqb{0,1}, 
\end{align}
where $N_x$ is the number of degrees of freedom in $x$-direction.  
The corresponding root mean square (RMS) profile of \eqref{eq:mean_velo}, which is a measure for the deviation from the mean profile, is given by 
\begin{align*}
	u_{1,\rms}\rb{t^*, y}  
		= \rb{\frac1{N_x}
			\sum_{i=1}^{N_x} \sqb{ 
				u_1\rb{t^*, x_i, y } - \bra{u_1}\rb{t^*, y}
			}^2
		}^{1/2}, \quad y\in\sqb{0,1}.
\end{align*}
Note that the degrees of freedom at the periodic boundary are used only once in the sums. 
Mean velocity profiles at certain points in time or even over the whole time interval, i.e.,  
\begin{align*}
	\bra{u_1}\rb{y}  = \frac{1}{N_t} \sum_{i=1}^{N_t} \bra{u_1}\rb{t_i, y} 	
\end{align*}
in the case of equidistant time steps, sometimes together with the associated RMS profile, were studied in  \cite{BoersmaEtAl97,GravemeierEtAl05,AhmedEtAl17}.
Here, $N_t$ denotes the number of time steps.
Due to the fact that in order to obtain the velocity profiles a lot of averaging has to be performed, it turns out that it
is not particularly demanding to achieve mesh convergence for this quantity.
This can already be seen in the literature.
Therefore, we will not plot velocity profiles again in this work.

%---------
\emph{Kinetic energy.}
The most frequently monitored quantity of interest is the kinetic energy of the 
flow, given by
\begin{align*}
	\Kin\rb{t^*,\uu} 
		= \frac{1}{2} \norm{\uu\rb{t^*}}_\LP{2}{\OMEGA}^2.
\end{align*}
For the studied flow problem, an energy inequality holds for the velocity field (which does not become constant) and therefore, the physically correct behavior of $\Kin$ is that it strongly monotonically decreases. 
If in practical computations of this problem, at some point in time, $\Kin$ increases again, there is something upsetting the energy balance of the scheme.
The kinetic energy was studied in \cite{VremanGeurtsKuerten97,GriebelKoster00,GravemeierEtAl05,John05,Burman07,AhmedEtAl17,SchroederLube18,YangBadiaCodina16}.
As an easily computable measure for the amount of energy in the flow, we include $\Kin$ in the set of  evaluated quantities of interest in this paper. 
In \cite{YangBadiaCodina16}, also the temporal change of the kinetic energy 
\begin{align*}
	-\frac{\drm \Kin\rb{t^*,\uu}}{\drm t}	
\end{align*}
was monitored. 

%---------
\emph{Kinetic energy spectra.}
The longitudinal Fourier transform of the streamwise velocity component $u_1$ is given by 
\begin{align*}
	\widehat{u}_1\rb{\kappa,t^*,y}
	= \int_0^1 u_1\rb{t^*,x,y} \exp\rb{-\im \kappa x} \dx,  
	\quad y\in\sqb{0,1},	
\end{align*}
where $\kappa \in \mathbb N$ is the wavenumber. Then, the longitudinal spectrum of $u_1$ is defined by 
\begin{align*}
	E(\kappa, t^*) 
		= \int_{0}^{1} \abs{ \widehat{u}_1 \rb{\kappa,t^*,y} }^2 \dy.
\end{align*}
Kinetic energy spectra at various points in time were studied in detail in \cite{LesieurEtAl88}; further results can be found in \cite{BoersmaEtAl97,SchneiderFarge00,IannelliEtAl03}. 
In this work, we present some longitudinal kinetic energy spectra for $\Rey=\num{10000}$.

%---------
\emph{Vorticity.}
As already mentioned, a fundamental feature of the Kelvin--Helmholtz instability problem is the development of vortices. Isolines or surface plots of the vorticity
\begin{align*}
	\omega 
		= \nabla\times \uu
		= \partial_x u_2 - \partial_y u_1	
\end{align*}
are an appropriate way for visualizing the flow field. 
Pictures presenting the vorticity can be found in every publication with simulations of the Kelvin--Helmholtz instability problem. 
For example, vorticity snapshots are shown in the qualitative study \cite{BurmanEtAl17}.
In Section \ref{sec:ComputationalStudies}, we provide exactly such plots.

%---------
\emph{Vorticity thickness.}

A popular quantity of interest that is monitored is the vorticity thickness defined by 
\begin{align} \label{eq:VortThick}
	\delta(t^*) 
	= \frac{2 u_\infty}{\sup_{y\in\sqb{0,1}} \abs{\bra{\omega}\rb{t^*,y}}}, 
\end{align}
where $\bra{\omega}$ is the integral mean vorticity in periodic direction, that is, 
\begin{align*}
	\bra{\omega}\rb{t^*,y} 
		= \frac{\int_{0}^1 \omega\rb{t^*,x,y}\dx }{\int_{0}^1 \dx} 
		=  \int_{0}^1 \omega\rb{t^*,x,y} \dx.
\end{align*}
It is common to monitor the vorticity thickness relative to the initial vorticity thickness $\delta_0$, i.e., $\delta(t^*)/\delta_0$.
Results for the evolution of the vorticity thickness over time are presented in 
\cite{LesieurEtAl88,NageleWittum03,GravemeierEtAl05,John05,YangBadiaCodina16,AhmedEtAl17,SchroederLube18} and also in this work.
Here, in order to compute the supremum in \eqref{eq:VortThick} numerically, we consider the maximum of $\bra{\omega}$ over 1024 equidistant horizontal lines in $\sqb{0,1}$.
In case that the evolution of the vorticity thickness shows an oscillatory behavior, the corresponding vortices in the Kelvin--Helmholtz problem are of ellipsoidal shape.
Viscosity effects act as a smoothing mechanism which is responsible for attenuating the oscillations in this quantity.
If the vorticity thickness remains more or less constant over time, the corresponding vortices are rather circular.

%---------
\emph{Enstrophy.} 
Mathematically, enstrophy is defined by 
\begin{align*}
	\Ens\rb{t^*,\uu} 
		= \frac{1}{2}\norm{ \nabla \times \uu\rb{t^*}}_\Lp{2}{\OMEGA}^2
		= \frac{1}{2} \norm{\omega\rb{t^*}}_\Lp{2}{\OMEGA}^2.
\end{align*}
Results with respect to the temporal evolution of the enstrophy for the Kelvin--Helmholtz instability problem are presented in \cite{SchneiderFarge00,GriebelKoster00,OnateEtAl07,SchroederLube18} and also in this work. 
Also for the enstrophy, the physically correct behavior is a monotone decrease from its initial value, as will be explained in more detail in 
Section~\ref{sec:SelfOrga2DTurbulence}.

%---------
\emph{Palinstrophy.}
Mathematically, palinstrophy is defined by 
\begin{align*}
	\Pal\rb{t^*,\uu} 
		= \frac{1}{2}\norm{ \nabla\rb{\nabla \times \uu\rb{t^*}}}_\LP{2}{\OMEGA}^2
		= \frac{1}{2} \norm{\nabla\omega\rb{t^*}}_\LP{2}{\OMEGA}^2.
\end{align*}
In the context of 2D turbulence and self-organization, see Section \ref{sec:SelfOrga2DTurbulence}, this is a very important quantity because it drives the dissipation process in 2D. 
Therefore, we will evaluate $\Pal$ for our reference results.
Note that $\Pal$ can indeed spontaneously increase in time; cf. \cite[Sec.~3.3]{DoeringGibbon95}. 
Also, it is important to mention that in any form of finite element analysis, this quantity is not controlled theoretically if the problem is discretized in velocity-pressure variables.
Thus, the evolution of the palinstrophy represents the most challenging quantity which is to be monitored.
We monitor also the palinstrophy with the aim of being able to assess the dynamics of the underlying problem more accurately.
Furthermore, recently, an increased interest, both theoretically and numerically, in using all the quantities $\Kin$, $\Ens$ and $\Pal$ for analyzing merging processes in incompressible 2D flows at high Reynolds numbers can be observed \cite{GarganoEtAl11,AyalaProtas14,AyalaProtas14b,ClercxVanHeijst17}.
This fact further underlines the importance of evaluating them.

%---------
\emph{Time intervals of the pairings.} 
In a Kelvin--Helmholtz instability problem with four primary eddies, a first pairing occurs where each two of them create a new larger eddy. 
Later, these two new eddies pair to a big final eddy. 
Important quantities of interest are the times resp.\ time intervals of these two pairings. 
We will see that especially the occurrence of the last pairing in time cannot be reliably predicted.
Besides from the visualization of the vorticity, the pairing times were determined in the literature by increases of the vorticity thickness. 
In addition, the pairings can be also observed by peaks of the palinstrophy and/or, at least for higher Reynolds numbers, a decrease in the enstrophy.
As it seems difficult to quantify the times or time intervals for the pairing process in a meaningful way, we rather rely on a comparison of the palinstrophy and enstrophy which indicate the occurrence of vortex pairings. 

%------------------------------------------------------------------------------------------------
%-----------------------------SELF ORGANIZATION AND 2D TURBULENCE--------------------------------
%------------------------------------------------------------------------------------------------
\section{Self-organization and 2D turbulence}
\label{sec:SelfOrga2DTurbulence}

There is a mathematical theory for the self-organization of solutions of the 2D Navier--Stokes equations developed by van Groesen in \cite{VanGroesen88}.
To the best of our knowledge, the mathematical CFD community seems to be only little aware of this theory. 
Thus, for the interested reader, this theory will be sketched in this section, thereby extending it to the boundary conditions present in the Kelvin--Helmholtz instability problem. 
The most important conclusions from this theory for this problem are summarized in Remark \ref{rem:VanGroesen1} and \ref{rem:VanGroesen2}.

Consider the parameterized domain $\Omega= \rb{0,\frac{\pi}{a}} \times \rb{0,\frac{\pi}{b}} \subset \R^2$ with outer unit normal $\nn$ (for $a=b=\pi$ we again obtain the above introduced domain). 
We search for periodic solutions in $x$-direction with a no-penetration condition in the normal direction (no flow across the boundary in $y$-direction) and a free-slip for the tangential direction (no viscous stress along the boundary in $x$-direction), see \cite{Gunzburger89,ErnGuermond04}, i.e.,
\begin{align} \label{eq:nopenetr}
  \uu \ip \nn = 0 \quad \text{and} \quad  \rb{- \nu \nabla \uu \ip \nn} \times \nn = \zero \quad \text{on} \quad 
  \rb{ 0,\frac{\pi}{a} } \times \set{0,\frac{\pi}{b}},
\end{align}
For the two-dimensional Kelvin--Helmholtz problem, the motion is completely driven by the (regularized)
initial condition \eqref{eq:ini_cond}.
Using standard arguments (taking the curl), for the vorticity $\omega = \nabla \times \uu$ (scalar-valued in 2D), problem \eqref{eq:TINS} transforms to 
\begin{subequations}\label{eq:vorticity}
\begin{empheq}[left=\empheqlbrace]{alignat=2}  
  	\partial_t \omega - \nu \Delta \omega + \uu \ip \nabla \omega &= 0, 
  		\qquad &&\text{in }\rsb{0,\tend} \times \Omega, \\
  	\omega\rb{0,\cdot} &= \omega_0\rb{\cdot } 
  		\qquad &&\text{in } \Omega.	
\end{empheq}	
\end{subequations}
Using $\DIV \uu = 0$ and \eqref{eq:nopenetr}, we now consider mixed boundary conditions
\begin{subequations} \label{eq:mixed-bv}
\begin{empheq}{alignat=2} 
	\omega &  = 0 \quad 
        \text{on}\quad
        \rb{ 0,\frac{\pi}{a} } \times \set{0,\frac{\pi}{b}}, \\
  	\omega\rb{0,y} & =  \omega\rb{\frac{\pi}{a},y}, \quad 
  		\rb{\nabla \omega \ip \nn} \rb{0,y} 
  			+\rb{\nabla \omega \ip \nn}\rb{\frac{\pi}{a},y} = 0
   			\quad  \text{on} \quad 0 \leqslant y \leqslant \frac{\pi}{b},
\end{empheq}	
\end{subequations}
where $\nn(0,y) = - \nn(\frac{\pi}{a},y),~y\in[0,\frac{\pi}{b}]$.
It is straightforward to show that the vorticity problem \eqref{eq:vorticity}--\eqref{eq:mixed-bv} admits a unique solution $\omega$. 

Using arguments from the theory of dynamical systems, Van Groesen presented in  \cite{VanGroesen88} a very convincing theory for the self-organization of solutions to the 2D Navier--Stokes equations. 
More precisely, he analyzed the vorticity problem \eqref{eq:vorticity} with homogeneous Dirichlet condition $\omega =0$ on $\partial \Omega$. 
Our goal here is to adapt the theory by Van Groesen to the incompressible Kelvin--Helmholtz problem \eqref{eq:vorticity} with mixed boundary conditions \eqref{eq:mixed-bv}. 

Recall the definitions of kinetic energy $\Kin$, enstrophy $\Ens$ and palinstrophy $\Pal$ (now with argument $\omega$ instead of $\uu$):
\begin{align*} 
	\Kin\rb{\omega} = \frac{1}{2} \int_\Omega \abs{\uu}^2 \dxx, \quad 
	\Ens\rb{\omega} = \frac{1}{2} \int_\Omega \abs{\omega}^2 \dxx, \quad 
	\Pal\rb{\omega} = \frac{1}{2} \int_\Omega \abs{\nabla \omega }^2 \dxx.
\end{align*}
Please note that $\Kin$ and $\Ens$ are invariants of the incompressible Euler problem, i.e.\ of \eqref{eq:TINS} with $\nu=0$.
Following \cite{VanGroesen88}, we consider critical (i.e.\ stationary) points of $\Ens$ on level sets of $\Kin$; that is,
\begin{align} \label{eq:crit}
	\crit \set{ \Ens\rb{\omega}\colon \Kin\rb{\omega} = \gamma; 
		~\text{with mixed BCs } \eqref{eq:mixed-bv};~\gamma > 0 \text{ fixed} } .
\end{align}
Since
$\Ens$ and $\Kin$ are quadratic functionals, the solutions of the constrained minimization problem \eqref{eq:crit} are solutions of \eqref{eq:crit} with $\gamma=1$ multiplied by the scaling factor $\sqrt{\gamma}$. 
Moreover, they are, up to scaling, the critical points
\begin{align*}
	\crit \set{ \Quo\rb{\omega}= \frac{\Ens\rb{\omega}}{\Kin\rb{\omega}}\colon 
   		~\text{with mixed BCs } \eqref{eq:mixed-bv};~ \gamma \not\equiv 0 } 	
\end{align*}
for the Rayleigh quotient $\Quo\rb{\omega}$.
These critical points are precisely the solutions of the eigenvalue problem
\begin{align}  \label{eq:EVP}
	-\Delta \omega = \lambda \omega 
  		\quad \text{in }\Omega, \quad  
  		\text{with mixed BCs } \eqref{eq:mixed-bv},
\end{align}	
see, for example, \cite[Section 44.5]{Zeidler85}.
Using the separation method, for $\kk = \rb{k_1,k_2} \in \N \times \N_0$ the solutions can be written as
\begin{align*}
	\lambda_\kk = 4ak_1^2+ bk_2^2, \quad 
	\widehat{w}_\kk^{(1)} = \frac{2\sqrt{ab}}{\pi} \sin\rb{2ak_1x} \sin\rb{bk_2y}, \quad 
	\widehat{w}_\kk^{(2)} = \frac{2\sqrt{ab}}{\pi} \cos\rb{2ak_1x} \sin\rb{bk_2y}.	
\end{align*}
The eigenfunctions $\widehat{w}_\kk^{(1,2)}$ form a complete $\Lp{2}{}$-orthonormal set with $\Kin\rb{\widehat{w}_\kk}=\rb{2\lambda_\kk}^{-1}$ and $\Ens\rb{\widehat{\omega}_\kk} = 1/2$. 
Denote the increasing set of eigenvalues by $\mu_k$ with $0 < \mu_1 = \lambda_{1,0} < \mu_2 < \mu_3 < \ldots $ and by $E_k$ the eigenspace corresponding to $\mu_k$. 
Note that the smallest eigenvalue $\mu_1= \lambda_{1,0}$ is simple.

\begin{figure}[t]
	\centering
	\includegraphics[width=0.9\textwidth]{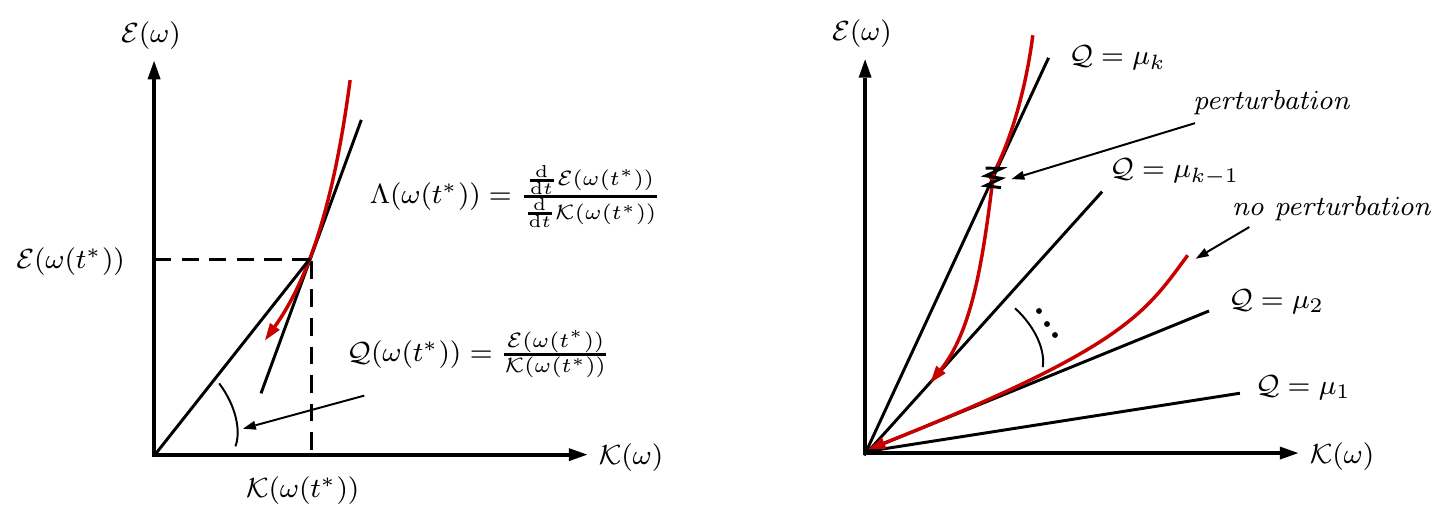}
	\caption{Sketch of the $\Kin$-$\Ens$-plane with Rayleigh quotient $\Quo$ and dissipation rate quotient $\Lambda$ (left). Different curves in the $\Kin$-$\Ens$-plane which correspond to the invariant sets and their asymptotic instability of Lemma~\ref{lem:InvariantSets} (right).}
	\label{fig:Self-organisation-sketches}
\end{figure}

As already mentioned, $\Kin$ and $\Ens$ are invariants of the motion for the incompressible Euler problem. 
More precisely, direct verification shows the following result.
%-----------------
\begin{thmLem}
Any $w \in E_k, k \in {\mathbb N}$, is a time-independent solution of the incompressible Euler problem \eqref{eq:TINS} with $\nu= 0$. 
Moreover, the so-called `Taylor vortices' $w \exp\rb{-\nu \mu_kt}$ are exact solutions of the incompressible Navier--Stokes problem \eqref{eq:TINS}.
\label{lem:EigFctSol}
\end{thmLem}
%-----------------
Following \cite{VanGroesen88}, we are now looking for the motion of solutions in the $\Kin$-$\Ens$-plane. 
The vorticity $\omega = \omega\rb{t}$ in \eqref{eq:vorticity} defines a continuous curve 
\begin{align} \label{eq:phasemotion}
	t \mapsto \sqb{\Kin\rb{\omega (t)},\Ens\rb{\omega (t)}}.
\end{align} 
Following \cite[Sec.~10.1.1]{Davidson04}, we obtain 
\begin{align} \label{eq:EW}
	\frac{\drm}{\drm t} \Kin\rb{\omega} = - 2 \nu \Ens\rb{\omega}, \qquad 
   	\frac{\drm}{\drm t} \Ens\rb{\omega} = - 2 \nu \Pal\rb{\omega} .
\end{align}
Moreover, this allows to introduce the so-called `dissipation rate quotient'
\begin{align*}
	\Lambda\rb{\omega} 
	= \frac{\Pal\rb{\omega}}{\Ens\rb{\omega}}
	= \frac{\frac{\drm}{\drm t} \Ens\rb{\omega}}{\frac{\drm}{\drm t} \Kin\rb{\omega}}. 	
\end{align*}
In the $\Kin$-$\Ens$ diagram, the quotients $\Quo$ and $\Lambda$ have a clear geometric interpretation; cf.\ Figure \ref{fig:Self-organisation-sketches} (left). 
At a fixed point of the curve \eqref{eq:phasemotion}, $\Quo$ is the angle of this point with the positive $\Kin$-axis whereas $\Lambda$  is the tangent to the curve \eqref{eq:phasemotion} at this point.
Furthermore, $\Lambda$ and $\Quo$ are homogeneous of degree $0$; hence with $\widehat{\omega} = \omega/\norm{\omega}_\Lp{2}{\OMEGA}$ one has $\Quo\rb{\omega} = \Quo(\widehat{\omega})$ and $\Lambda\rb{\omega} = \Lambda\rb{\widehat{\omega}}$ for all $\omega$. 

%-----------------
%-----------------

Let us further remark that the palinstrophy plays a crucial role in the high wavenumber scaling theory for two-dimensional turbulence. 
Indeed, using a length scale $L$, by introducing the enstrophy dissipation rate $\chi=2\nu L^{-2} \Pal\rb{\omega}$, from dimensional analysis it follows that for high wavenumbers $\kappa$, the kinetic energy spectrum behaves according to $E\rb{\kappa}\sim\chi^{\nf{2}{3}}\kappa^{-3}$; cf.\ \cite{DoeringGibbon95}. 
An immediate consequence from the dynamical system \eqref{eq:EW} is that both kinetic energy $\Kin$ and enstrophy $\Ens$ decrease monotonically (exponentially) as $t \to \infty$. 
Indeed, the min-max theorem (or variational theorem) for eigenvalues leads to $\Pal\rb{\omega} \geqslant \mu_1 \Ens\rb{\omega}$ for all $\omega$. This implies $\frac{\drm}{\drm t} \Ens = - 2 \nu \Pal \leqslant -2 \nu \mu_1 \Ens$ and hence
\begin{align*}
	\Ens\rb{\omega (t)} 
		\leqslant \Ens\rb{\omega(0,\cdot)} \exp\rb{-2 \nu \mu_1 t}
		=\Ens\rb{\omega_0} \exp\rb{-2 \nu \mu_1 t}.
\end{align*}
Similarly $\Ens\rb{\omega} \geqslant \mu_1 \Kin\rb{\omega}$ implies 
\begin{align*}
	\Kin\rb{\omega (t)} 
		\le \Kin\rb{\omega(0,\cdot)} \exp\rb{-2 \nu \mu_1 t}
		= \Kin\rb{\omega_0} \exp\rb{-2 \nu \mu_1 t}.
\end{align*}

Consider now the behavior of the Rayleigh quotient $\Quo=\Ens/\Kin$ along a solution. 
Due to the restriction to 2D, the change of $\Quo$ is driven exclusively by the viscosity $\nu$ (no vortex stretching!). 
It turns out that $\Quo$ decreases monotonically for $t \to \infty$. 
One obtains 
\begin{align} \label{eq:Q-punkt}
	\frac{\drm}{\drm t}\Quo 
		= \frac{\drm}{\drm t} \frac{\Ens}{\Kin} 
		= \frac{-2 \nu \Pal\Kin - \Ens\rb{-2 \nu \Ens}}{\Kin^2} 
		= - 2 \nu \rb{ \frac{\Pal}{\Kin} - \frac{\Ens^2}{\Kin^2} } 
		= - 2\nu \Quo (\Lambda - \Quo) .
\end{align}

%-----------------
\begin{thmLem}
For each solution $\omega$ of \eqref{eq:vorticity} with mixed boundary conditions \eqref{eq:mixed-bv}, one has $\Lambda\rb{\omega}-\Quo\rb{\omega} \geqslant 0$. 
Moreover, $\Lambda\rb{\omega} - \Quo\rb{\omega} = 0$ is valid iff $\omega \in E_k$ for some $k \in \N$.
\label{lem:LambdaQ}
\end{thmLem}
%-----------------
\begin{thmProof}
This follows as in the proof of Lemma 6.1 in \cite{VanGroesen88}.
\end{thmProof}
%-----------------

%-----------------
\begin{thmLem}
For each solution $\omega(t)$ of \eqref{eq:vorticity} with initial condition $\omega_0$, one has $\frac{\drm}{\drm t} \Quo\rb{\omega (t)} \leqslant 0$.
Moreover, the limit 
\begin{align*}
	q\rb{\omega_0} = \lim_{t \to \infty} \Quo\rb{\omega(t)}
\end{align*} 
exists.
This limit satisfies $\mu_1 \leqslant q\rb{\omega_0} \leqslant \Quo\rb{\omega_0}$ with $q\rb{\omega_0} = \Quo\rb{\omega_0}$ iff $\omega_0 \in E_k$ for some $k \in \N$ and this equality holds iff $\omega\rb{t}$ is a planar Taylor vortex.
Finally, one obtains
\begin{align*}
   q\rb{\omega_0} = \lim_{t \to \infty} \Lambda\rb{\omega (t)}.
\end{align*}
\end{thmLem}
%-----------------
%-----------------
\begin{thmProof} (Sketch)
We follow the proof of Corollary~6.2 in \cite{VanGroesen88}.
Since $\Quo \geqslant \mu_1$, from \eqref{eq:Q-punkt} and Lemma~\ref{lem:LambdaQ}, it follows that $\Quo$ decreases monotonically for any solution that is not a planar Taylor vortex.
Then $\Quo\rb{\omega(t)}$ has a limit as $t \to \infty$, since $\Quo$ is bounded from below.
Because $\Quo$ decreases monotonically, one has $\lim_{t \to \infty} \frac{\drm}{\drm t}\Quo = 0$, hence $\lim_{t \to \infty} \sqb{\Lambda\rb{\omega (t)} - \Quo\rb{\omega(t)} } =0$. 
\end{thmProof}
%-----------------

%-----------------
\begin{thmLem} \label{lem:IdealSituation}%
For each initial condition $\omega_0$ there is some $k \in \N$ such that for the limit of the Rayleigh quotient  $q\rb{\omega_0}= \mu_k$ holds for one of the eigenvalues $\mu_k$ of \eqref{eq:EVP}.	
\end{thmLem}
%-----------------
\begin{thmProof}
We refer to the technical proof of Proposition~6.3 in \cite{VanGroesen88} (see also Sections~7 and~9).
\end{thmProof}
%-----------------

We are now looking for sets which remain invariant for the flow. 
For $k \in \N$, let us define
\begin{align*}
	I_k = \set{ \omega \colon q\rb{\omega} = \lim_{t \to \infty} \Quo\rb{\omega(t)}= \mu_k }.
\end{align*}
Any set $I_k$ is invariant and contains the eigenspace $E_k$. 
Moreover, the intersection with any $I_j, j \ne k$, is empty while the union of all $I_k$'s is the whole function space. 
The following result characterizes the set $I_1$ and shows that $I_1$ is the only \emph{stable} one.
%-----------------
\begin{thmLem} \label{lem:InvariantSets}%
The invariant set $I_1$ contains the sets 
\begin{align*}
	\set{ \omega \colon \Quo\rb{\omega} < \mu_2 } 
  	\quad \text{and} \quad
  	\set{ \omega \colon \Quo\rb{\omega} = \mu_2,~ \omega \not\in E_2  } .	
\end{align*}
Moreover, each of the sets $I_k$ with $k \ne 1$ is asymptotically unstable in the following sense:
If $\omega\rb{t} \in I_k, k \ne 1$, then for any arbitrary small perturbation $\xi \in E_1 \cup E_2 \cup  \ldots \cup E_{k-1}$ there exists a point in time $\tilde{T}>0$ (sufficiently large) such that the solution with initial condition $\omega\rb{\tilde{T}} + \xi$ belongs to $I_j$ for some $j <k$.
\end{thmLem}
%-----------------

Note that the requirement $\xi \in E_1 \cup E_2 \cup \ldots \cup E_{k-1}$ means that the perturbation $\xi$ belongs to a coarser eigenspace than $E_k$, although $\norm{\xi}_\Lp{2}{\OMEGA}$ can be arbitrarily small.

%-----------------
\begin{thmProof}[of Lemma~\ref{lem:InvariantSets}]
See the proof of Proposition~6.4 in \cite{VanGroesen88}.	
\end{thmProof}
%-----------------

One can visualize the results in the $\Kin$-$\Ens$ plane; cf.\ Figure \ref{fig:Self-organisation-sketches} (right).
Each curve \eqref{eq:phasemotion} approaches some line $\Quo=\mu_k$ from above and is tangent to it at the origin. 
As a result of a small perturbation, such limiting line (with $k \ne 1$) can be crossed (in the absence of further perturbations). 
The decrease of $\Quo$ forces an asymptotic approach to some lower line $\Quo=\mu_l$ with $l<k$. 

In the perfectly unperturbed setting, Lemma \ref{lem:IdealSituation} predicts an idealized limit in which coherent structures (a certain number of vortices in the Kelvin--Helmholtz problem) form.
Lemma~\ref{lem:InvariantSets}, on the other hand, reacts according to nature and gives information about the structure of a perturbation and how strong it has to be in order to result in a different, perturbed limit.

As already pointed out at the very beginning of Section \ref{sec:Introduction}, nature is full of perturbations which inevitably lead to a behavior according to Lemma~\ref{lem:InvariantSets}.
Note that the initial vorticity to the initial data in \eqref{eq:ini_cond} is not $\Lp{2}{}$-orthonormal to $E_1$ so that even without any further perturbation the long time behavior will be dominated by the smallest eigenvalue, i.e.\ $\lim_{t \to \infty} \Quo\rb{\omega(t)}= \mu_1$. 

%-----------------
\begin{thmRem} \label{rem:VanGroesen1}
Note that Lemma~\ref{lem:InvariantSets} shows that there is only one stable state, but it does not contain any more information about it.
With respect to the Kelvin--Helmholtz problem, this stable set might correspond to the last vortex which forms at the end.
As already mentioned above, there is no consensus in the literature about the position of this last vortex as different discrete settings generally lead to different final states.
Let us emphasize that unfortunately, the theory presented here also does not make any prediction about the position of the last vortex.
Even more significantly, the instance in time where the last merging process occurs is, due to the sensitivity of the problem, not clear.
\end{thmRem}
%-----------------

%-----------------
\begin{thmRem} \label{rem:VanGroesen2}
The presented theory on the continuous level displays how sensitive flow simulations in 2D with high Reynolds numbers are.
One can expect that this property will be passed on to any attempt of obtaining discrete approximations to such systems.
Indeed, in every numerical method there are inevitably errors and perturbations.
And as we will show in Sections~\ref{sec:ComputationalStudies} and~\ref{sec:Perturbations}, those perturbations coming from discretizations can have a very dramatic effect on numerical approximations of the \emph{evolution} of the Kelvin--Helmholtz instability problem.
\end{thmRem}
%-----------------

%------------------------------------------------------------------------------------------------
%---------------------------------------DIV-FREE HDG---------------------------------------------
%------------------------------------------------------------------------------------------------
\section{High-order divergence-free IMEX HDG methods}
\label{sec:DivFreeHDG}

In this section, we briefly want to comment on the discretization, both in space and in time, we used to obtain our numerical solutions.
As this is not intended to be the focus of the present paper, we try to only explain the crucial ingredients and philosophies involved.
The discretization is based on the primitive velocity/pressure variables \eqref{eq:TINS} and not, for example, on the vorticity equation or related coupled problems.
All computations in this work have been carried out using the high-order finite element library \texttt{NGSolve} \cite{Schoeberl14}. 
The computations have been done on several different computers and, in doing so, we carefully verified and guaranteed the reproducibility and consistency between them. 

%-----------------------------------SPACE DISCRETIZATION-----------------------------------------
\subsection{Space discretization}

\begin{figure}[t]
	\centering
	\includegraphics[width=0.8\textwidth]{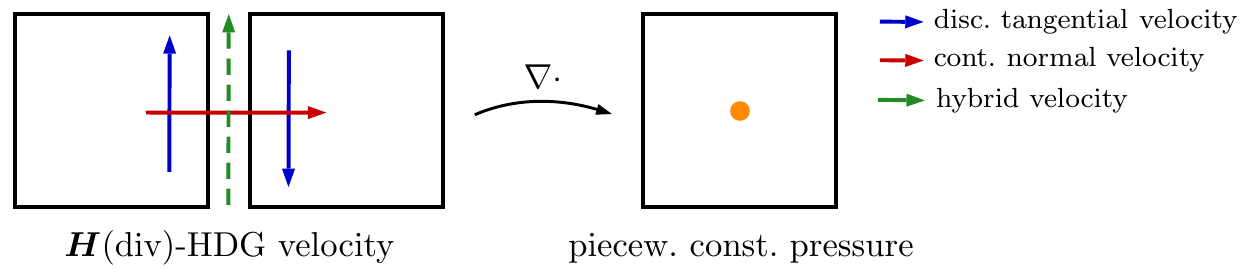}
	\caption{Sketch of the used HDG velocity/pressure pair on squares.}
	\label{fig:HDG-sketch}
\end{figure}

The space discretization which is used for the simulations in Section~\ref{sec:ComputationalStudies} is based on \cite{LehrenfeldSchoeberl16}.
Thus, we use a high-order, exactly divergence-free hybrid discontinuous Galerkin (HDG) method based on $\HDIV$ finite elements. 
The discretization of the viscous term is based on a hybrid version of the well-known symmetric interior penalty method.
We use `projected jumps', cf. Section~2.2.1 in \cite{LehrenfeldSchoeberl16}, static condensation and remove all higher-order non-divergence-free velocities from the basis such that the corresponding inf-sup stable discrete pressure space consists only of one constant per cell, cf.\ Remark~1 in \cite{LehrenfeldSchoeberl16}.
Figure~\ref{fig:HDG-sketch} displays a sketch of the resulting velocity/pressure pair.
Let us further mention that the chosen discretization is pressure-robust, which means that the velocity error is completely independent of the pressure error; cf.\ Remark~5 in \cite{LehrenfeldSchoeberl16}.
We choose this method as it -- compared to other standard discretization approaches -- combines features such as high order accuracy, important global and local conservation properties, energy-stability, polynomial, pressure and $\Rey$-semi-robustness, a minimal amount of numerical dissipation and computational efficiency \cite{LehrenfeldSchoeberl16,LedererSchoeberl2018,SchroederEtAl18}. 
Especially the combination of the robustness properties is hardly seen in other numerical discretization schemes.

%-----------------
%-----------------

Our numerical solutions in Section~\ref{sec:ComputationalStudies} are computed on meshes consisting of rectangles and therefore, the $\Hdiv$-conforming finite element space is based on the Raviart--Thomas element; cf.\ \cite{BoffiEtAl13}.
For all results, we use a higher-order approximation with $k=8$, where $k$ denotes the polynomial order that is completely included in the discrete velocity space.
We thus abbreviate the method by RT8.
In Section~\ref{sec:Perturbations}, on the other hand, results for computations on triangles are also shown.
On triangles, we use the Brezzi--Douglas--Marini element with $k=8$ (BDM8).

In order to deal with high Reynolds number flows, it is important to comment on the choice of numerical fluxes for the convection term.
As the finite element scheme is not based on $\HM{1}{}$-conforming discrete velocity spaces, the discrete velocity has discontinuities across interior inter-element boundaries.
In case of $\HDIV$ methods, only the tangential component of the velocity is discontinuous across those facets.
In the weak form, the \emph{upwinding} convection trilinear form we use then reads as follow:
\begin{align*} 
	c_h\rb{\uu_h;\uu_h,\vv_h}	
		=& \int_\Omega \rb{\uu_h\ip\nabla_h}\uu_h\ip\vv_h \dx
		- \sum_{F\in\Fi} \oint_F 
			\rb{\uu_h\ip\nn_F}\jmp{\uu_h}\ip\avg{\vv_h}  
			\ds \\
		&+ \sum_{F\in\Fi} \oint_F
			 \frac{1}{2} \abs{\uu_h\ip\nn_F}\jmp{\uu_h}\ip\jmp{\vv_h} 
			\ds, 			
\end{align*}	
where $\Fi$ is the set of interior facets, $\nn_F$ is a unit normal to a facet $F$, $\jmp{\cdot}$ and $\avg{\cdot}$ are the usual jump and average operators across facets and $\nabla_h$ is the broken (i.e.\ element-wise) gradient.
Naturally, facets belonging to a periodic edge of of $\partial\Omega$ are also collected in $\Fi$.
We notice that for piecewise polynomials of order $p$, in general, the numerically correct treatment of the convection form requires numerical quadrature of order $3p$, cf. Section~\ref{sec:InaccurateSolution}.
Especially, due to using Raviart--Thomas elements, we have $p=k+1$ for the volume integral.  
The last facet term incorporates upwinding in the DG sense; cf.\ \cite{PietroErn12,SchroederLube18}.
Note that whenever one deals with $\HDIV$ methods, jump and average operator actually only act on the tangential components. 
Note that the upwinding does not have any effect on sufficiently fine meshes.

For the computation of vorticity, enstrophy and palinstrophy we use broken (i.e.\ element-wise) derivatives; e.g.\ the enstrophy $\Ens\rb{t^*,\uu_h}$ is computed as
$\frac{1}{2} \norm{ \nabla_h \times \uu_h\rb{t^*}}_\Lp{2}{\OMEGA}^2 = \frac{1}{2} \sum_{T \in \T} \norm{ \nabla \times \uu_h\rb{t^*}}_\Lp{2}{\rb{T}}^2$ where $\T$ denotes the corresponding decomposition of the domain.

%-----------------
\begin{thmRem}
	In the literature \cite{GravemeierEtAl05,Burman07,AhmedEtAl17,YangBadiaCodina16} concerning Kelvin--Helmholtz simulations, due to different reasons, there are always stabilization mechanisms applied.
	Those stabilization mechanisms can be interpreted as additional numerical viscosity. 	
	As became evident in Section \ref{sec:SelfOrga2DTurbulence}, any small perturbation can lead to quantitative and qualitative differences in the computation of two-dimensional flows with high Reynolds numbers which is, in our opinion, the reason for decisively different solutions to Kelvin--Helmholtz instabilities in the literature.
	Especially, the merging of the last two vortices occurs relatively early in the existing literature.
	We believe that this is caused partially by the use of heavily stabilized methods.
	In this sense, the div-free HDG method we use here adds only a minimal amount of numerical dissipation from which our results clearly benefit \cite{SchroederEtAl18}.
\end{thmRem}
%-----------------

%------------------------------------TIME DISCRETIZATION-----------------------------------------
\subsection{Time discretization}
\label{sec:TimeDiscretization}

In this work, we decided to use an implicit-explicit (IMEX) time-stepping scheme where the nonlinear convection part is treated explicitly in time and the linear, stiff Stokes part is treated implicitly.
Thus, in each time step, only symmetric linear systems have to be solved, which can be done efficiently.
This efficiency makes it possible to choose a very small and constant time step size $\Delta t$ which, thereby, leads to accurate results as long as a resulting CFL condition is fulfilled.
Namely, in all the following computations $\Delta t=\delta_0\times\num{e-3}\approx \num{3.6e-5}$ has been chosen.
Note that the explicit treatment of the convection in time combines naturally with the standard DG (upwinding) formulation in space while the (`projected jumps') $\HDIV$-HDG formulation in space is specifically tailored to alleviate the computational costs of solving linear systems of Stokes-type operators \cite{LehrenfeldSchoeberl16}.
Let us remark that this certainly is only one possible approach for the temporal discretization.
Indeed, the numerical results in the next section should not depend on the particular time-stepping method used.

Any space discretization of the incompressible Navier--Stokes model leads to a nonlinear ODE system of the following form:
\begin{align}\label{eq:ODESystem}
	M\frac{\partial u}{\partial t}+Au+C\rb{u}u = F,
	\quad u\rb{0} = u^0.
\end{align}
Here, $M$ is a mass matrix, $A$ represents the Stokes bilinear form (viscosity and pressure), $C$ describes the nonlinear convection and $F$ stands for right-hand side forcing terms. For the present Kelvin--Helmholtz problem, the motion is completely determined by the initial condition $u^0$ and we have no forcing terms, i.e.\ $F=0$.

We chose the second-order semi-implicit BDF (SBDF2) method from \cite{AscherEtAl95} which combines a second-order BDF scheme with a second-order accurate extrapolation in time.
Supposing that $u^{n-1}$ and $u^n$ are already known, applying this time discretization to \eqref{eq:ODESystem} leads to solving the linear system
\begin{align*}
	\frac{1}{2\Delta t}M\sqb{3u^{n+1}-4u^n+u^{n-1}}
		= -2 C\rb{u^n}u^n + C\rb{u^{n-1}}u^{n-1} - Au^{n+1}.
\end{align*}
As usual for BDF methods, for the first time step a first-order implicit-explicit Euler method is used.
Rewriting the time-stepping scheme in incremental form results in the following approach for $n\geqslant 1$:
\begin{align*}
	\rb{M + \Delta t A}\sqb{u^1 - u^0}
		&= - \sqb{\Delta t C\rb{u^0}u^0 + \Delta t A u^0} ,\\
	\rb{M + \frac{2}{3}\Delta t A}\sqb{u^{n+1} - u^n}
		&= - \sqb{	
			\frac{4}{3}\Delta t C\rb{u^n}u^n -\frac{2}{3}\Delta t C\rb{u^{n-1}}u^{n-1} 
			+ \frac{2}{3}\Delta t A u^n - \frac{1}{3} M \rb{u^n-u^{n-1}}
		}.
\end{align*}
The abbreviation $M^\ast=M + \frac{2}{3}\Delta t A$ for the system matrix of the SBDF2 method is used.
We thus have to solve only one symmetric linear system per time step and the computational cost for this is determined by the structure of $M^\ast$.
This is the reason why we chose to use the SBDF class instead of, for example, the IMEX Runge--Kutta class of \cite{AscherEtAl97}.

Linear systems are solved with a sparse direct solver (sparse Cholesky from \texttt{NGSolve} \cite{Schoeberl14}) and iterative refinement.
For the chosen (relative) tolerance of $\num{e-12}$ (measured in the 2-norm) usually 2 or 3 iterative refinements are sufficient.
The gain in accuracy due to the iterative refinement proved to be very important for the Kelvin--Helmholtz instability problem; see also Section~\ref{sec:InaccurateSolution}.

%-----------------
\begin{thmRem}
  Concerning stability, it is worth mentioning that using IMEX schemes may lead to time step restrictions. 
In the context of this study, it turned out that the use of the considered IMEX method automatically enforces a sufficient temporal resolution once stability is provided, i.e.\ the time step restriction due to stability considerations and the time step restriction due to accuracy demands are in the same order of magnitude.
Choosing a time step smaller than $\Delta t=\delta_0\times\num{e-3}\approx \num{3.6e-5}$ only had a marginal impact on our results.
\end{thmRem}
%-----------------

%------------------------------------------------------------------------------------------------
%------------------------------------COMPUTATIONAL STUDIES---------------------------------------
%------------------------------------------------------------------------------------------------
\section{Computational studies}
\label{sec:ComputationalStudies}

In order to show how sensitive the Kelvin--Helmholtz instability problem is, we present results on a sequence of square meshes with $16^2$ to $256^2$ elements.
A feeling for the computational cost can be obtained in Table~\ref{tab:DOFs} where the resulting numbers of degrees of freedom (DOFs) and non-zero entries in the system matrix $M^\ast$ are summarized.

%-----------------
\begin{table}[b]
\caption{Overview of meshes, DOFs and non-zero entries of $M^\ast$ based on a discretization with RT8. The DOFs are counted before static condensation whereas the non-zero entries are counted from the system matrix after static condensation. The hybrid facet DOFs are not counted.}
\label{tab:DOFs}
\centering 
\begin{tabular}{crrrrr} 
\toprule
	Mesh	
		& $16^2$ 
		& $32^2$		
		& $64^2$	
		& $128^2$
		& $256^2$
		\\ 
\otoprule
	$\#\set{\uu\,\mathrm{DOFs}}$		
		& \num{21280} 	
		& \num{84544}
		& \num{337024} 
		& \num{1345792} 
		& \num{5378560} 
		\\
	$\#\set{p\,\mathrm{DOFs}}$		
		& \num{256} 	
		& \num{1024}
		& \num{4096} 
		& \num{16384} 
		& \num{65536} 
		\\		
	$\#\set{\text{nz}(M^\ast)}$		
		& \num{1075472} 	
		& \num{4292640}
		& \num{17152064} 
		& \num{68571264} 
		& \num{274211072} 
		\\
\bottomrule		
\end{tabular}
\end{table}
%-----------------

\begin{figure}[t]
\centering
\subfigure{\includegraphics[width=0.31\textwidth]{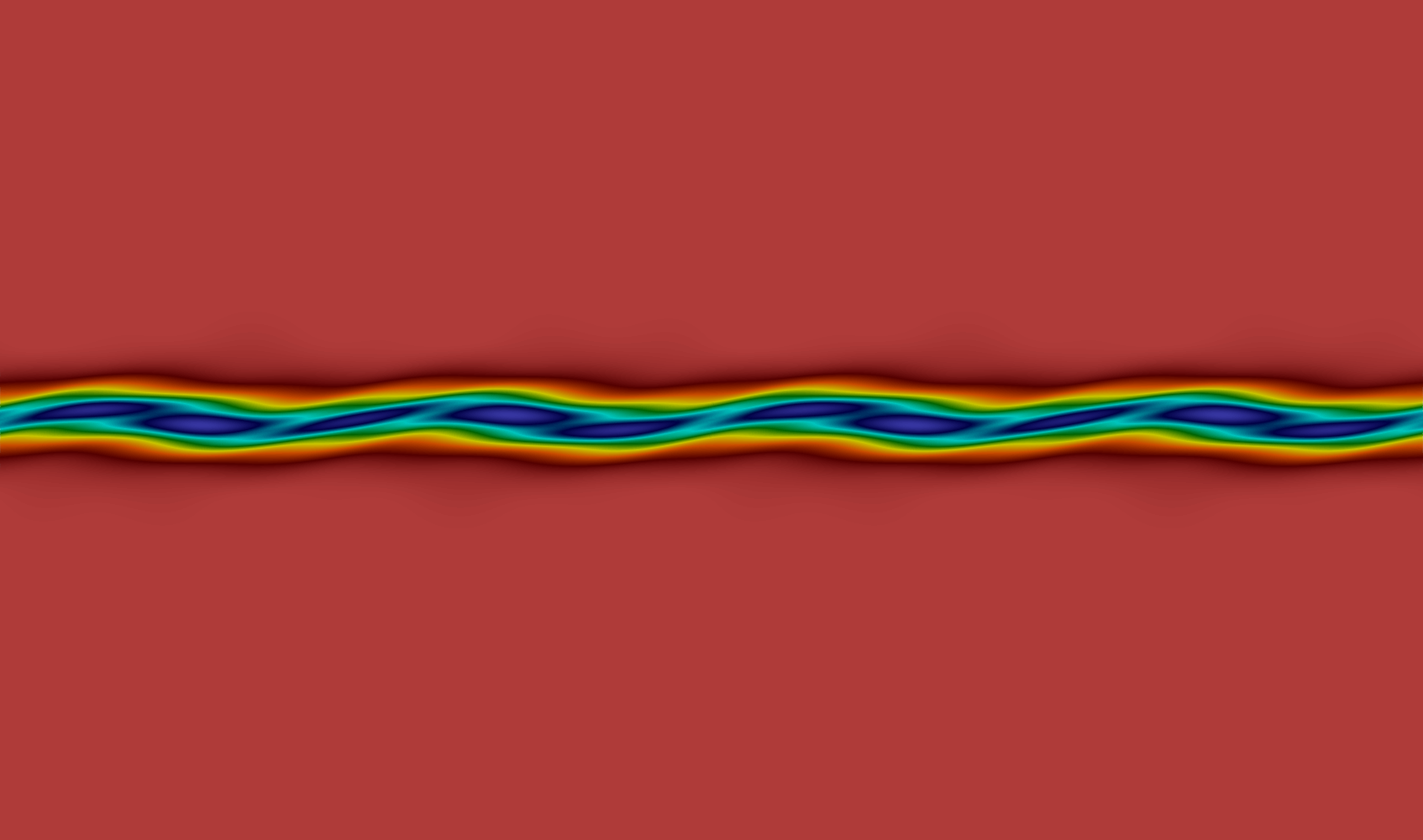}} \goodgap
\subfigure{\includegraphics[width=0.31\textwidth]{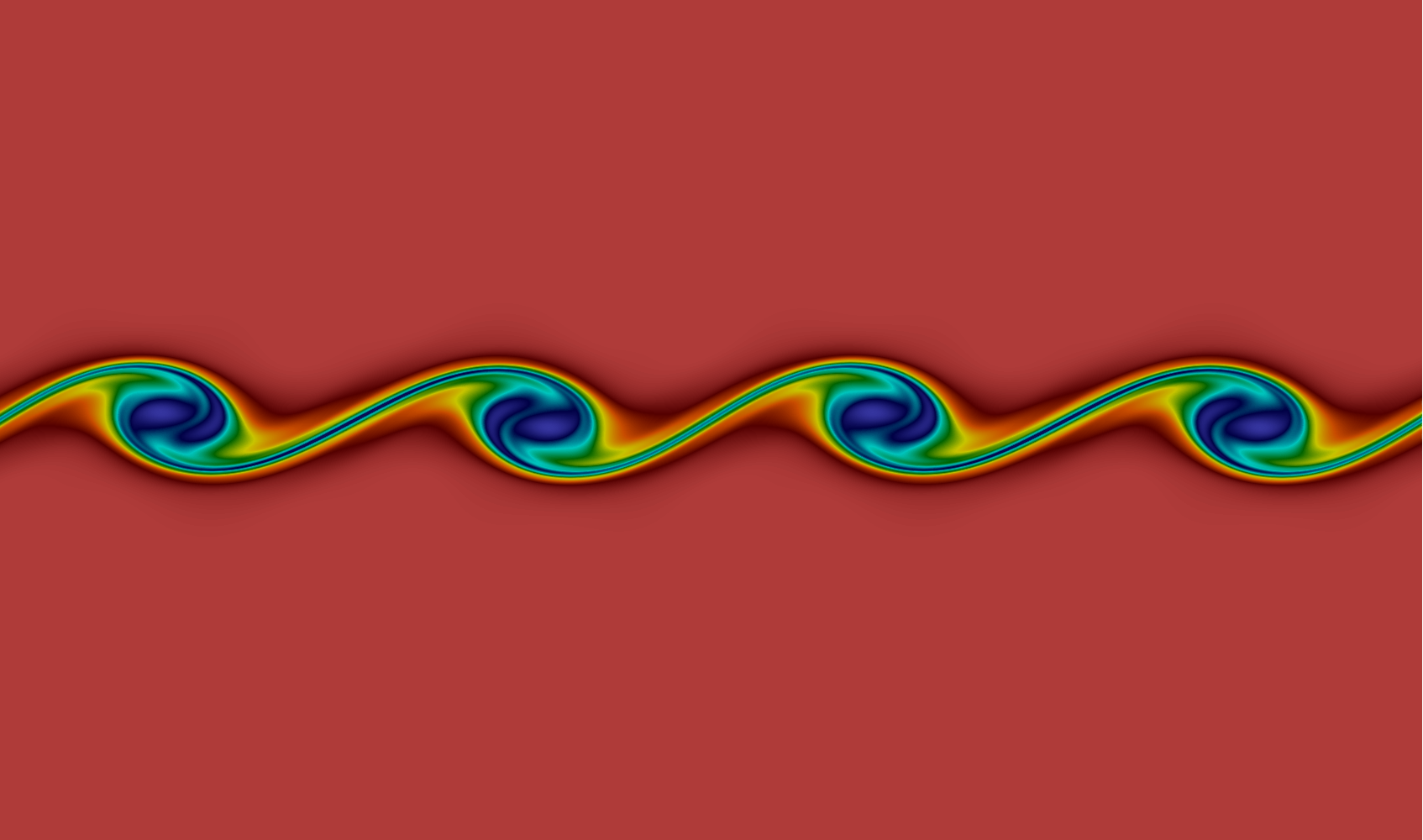}} \goodgap
\subfigure{\includegraphics[width=0.31\textwidth]{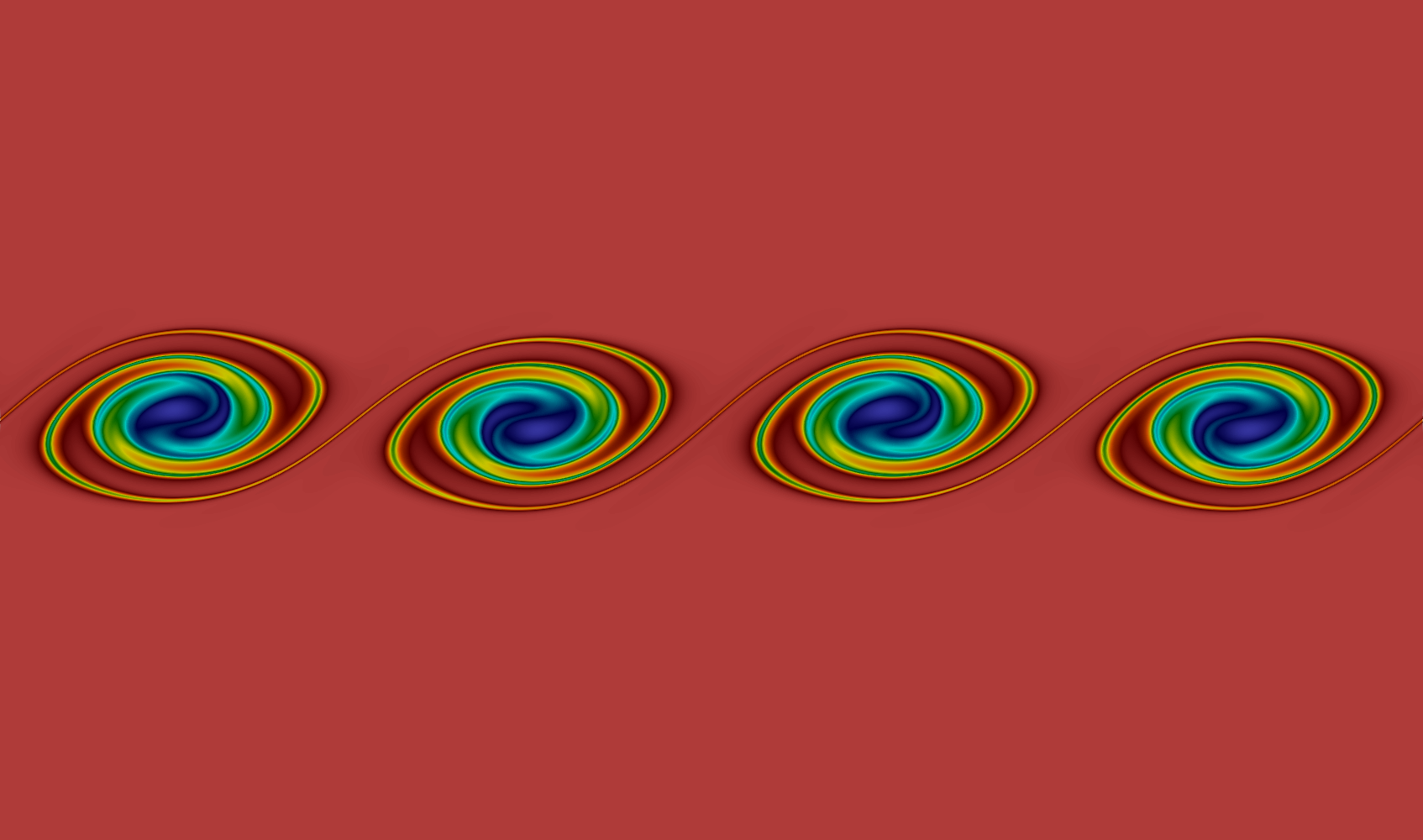}} 
\subfigure{\includegraphics[width=0.31\textwidth]{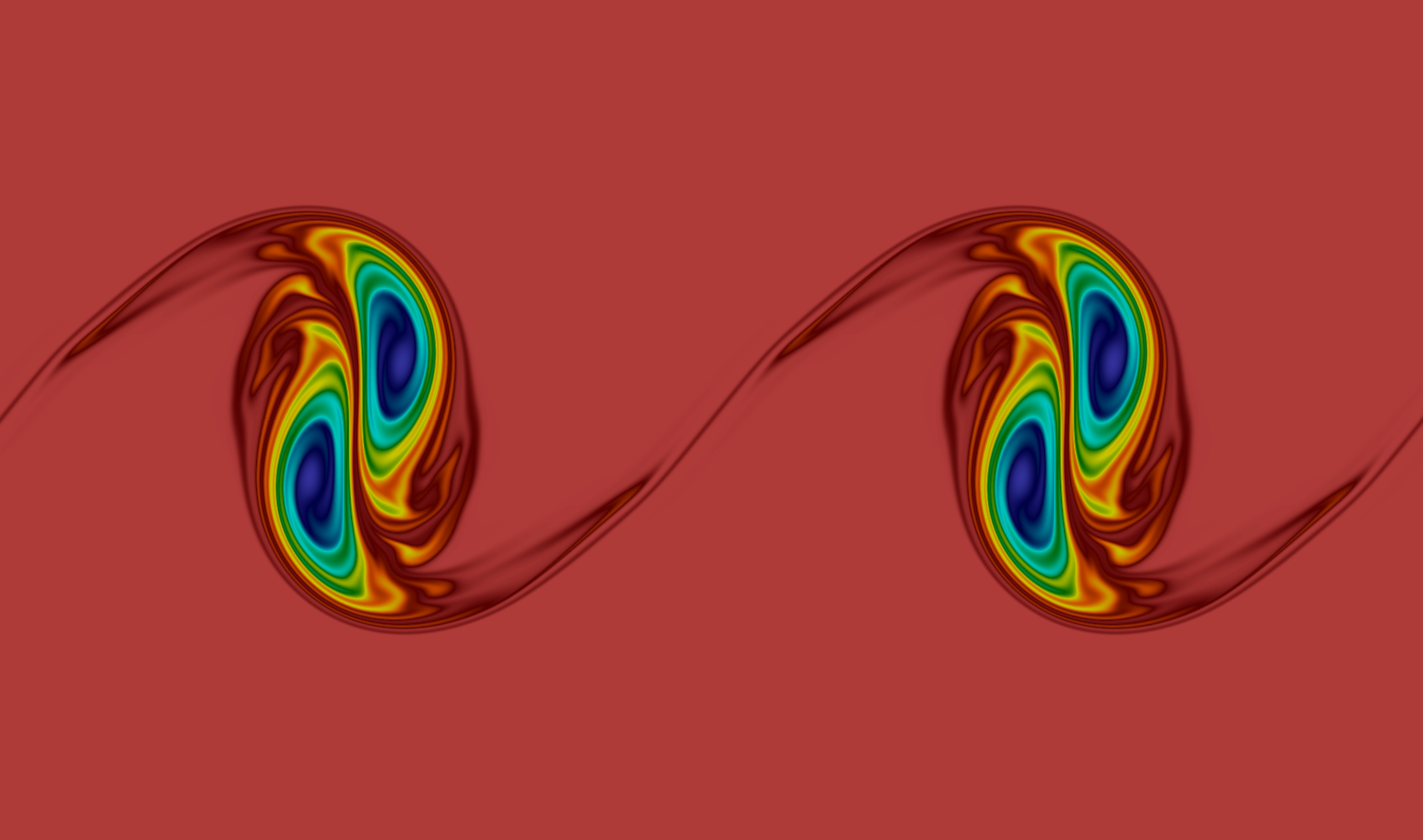}} \goodgap
\subfigure{\includegraphics[width=0.31\textwidth]{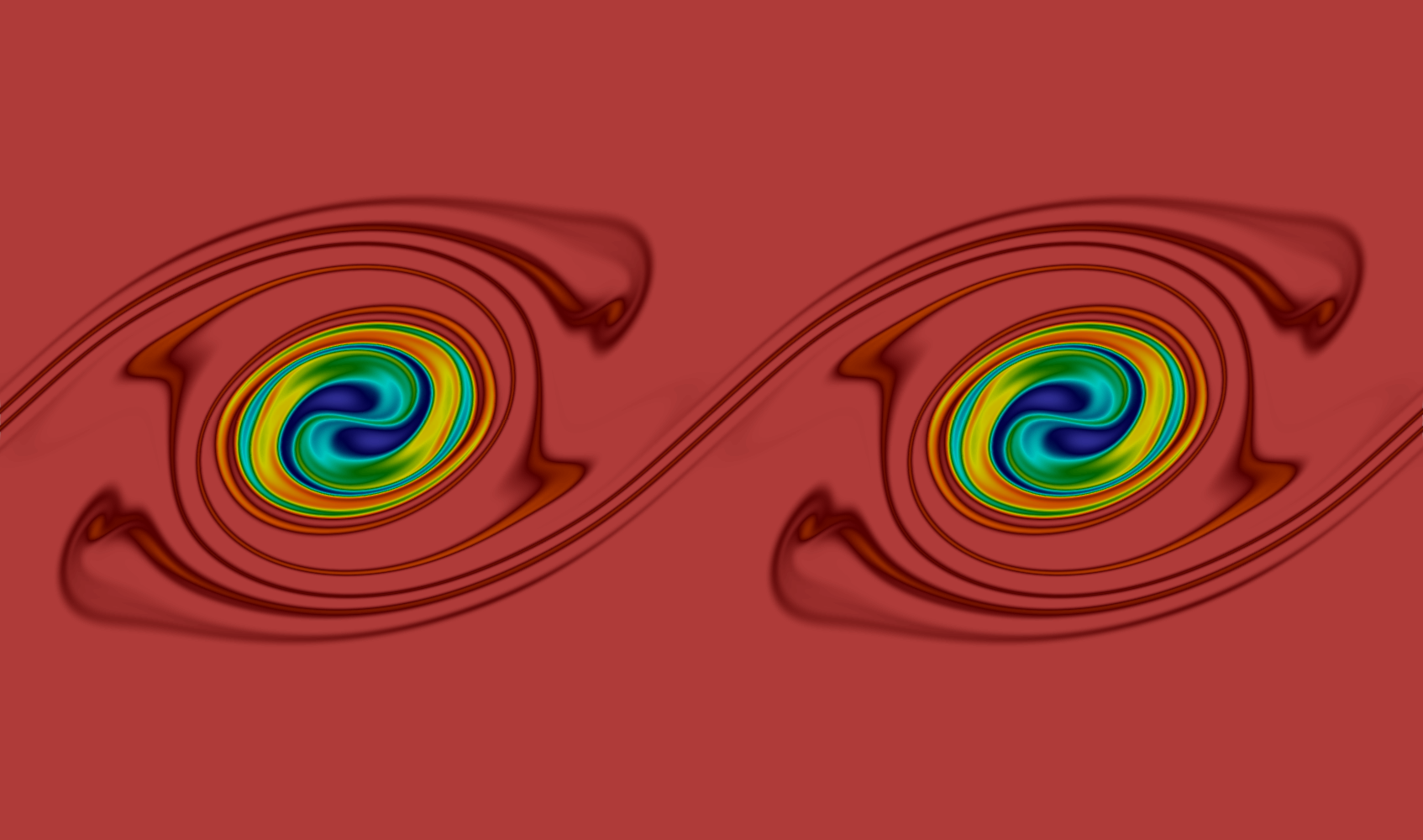}} \goodgap
\subfigure{\includegraphics[width=0.31\textwidth]{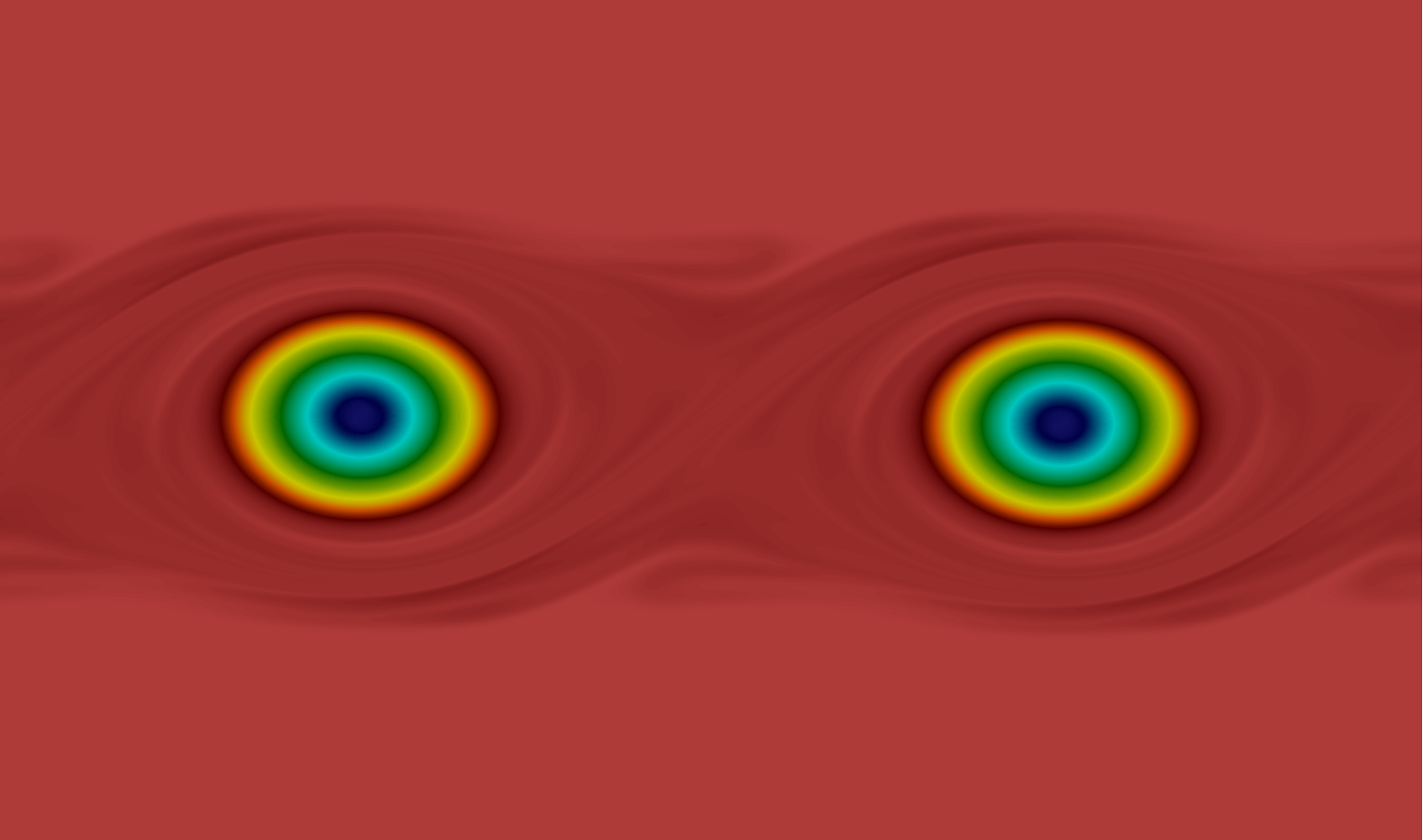}}
\subfigure{\includegraphics[width=0.31\textwidth]{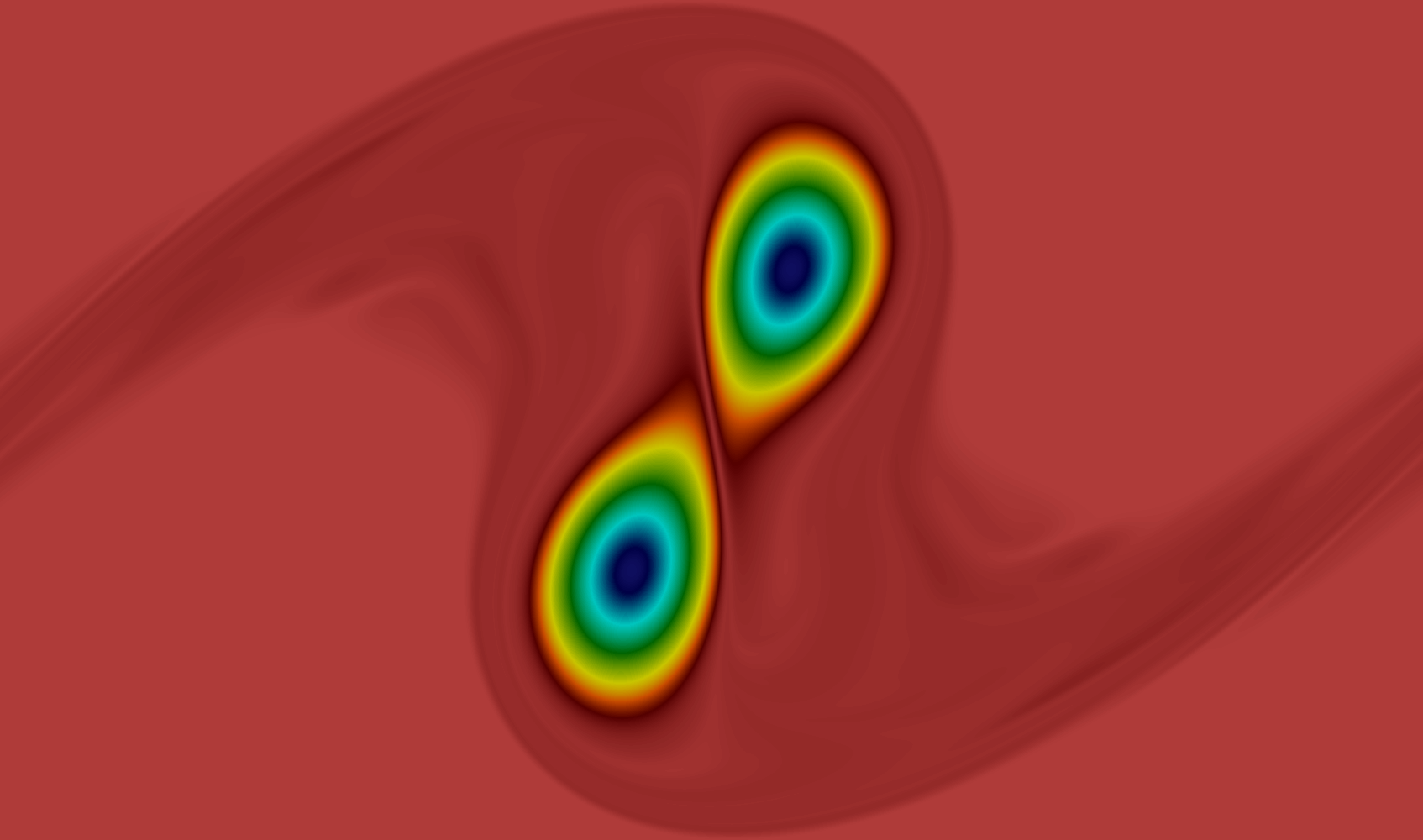}} \goodgap
\subfigure{\includegraphics[width=0.31\textwidth]{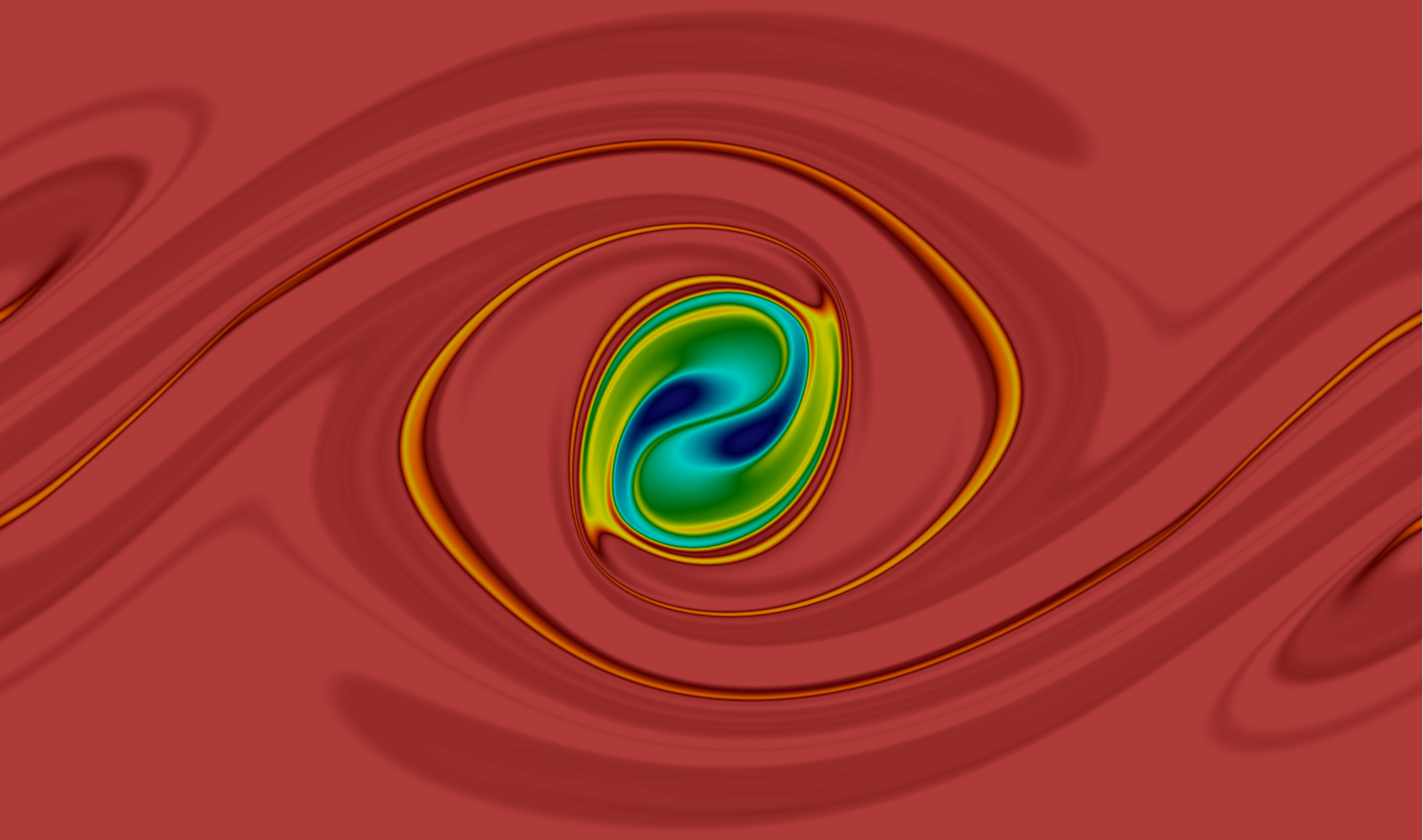}} \goodgap
\subfigure{\includegraphics[width=0.31\textwidth]{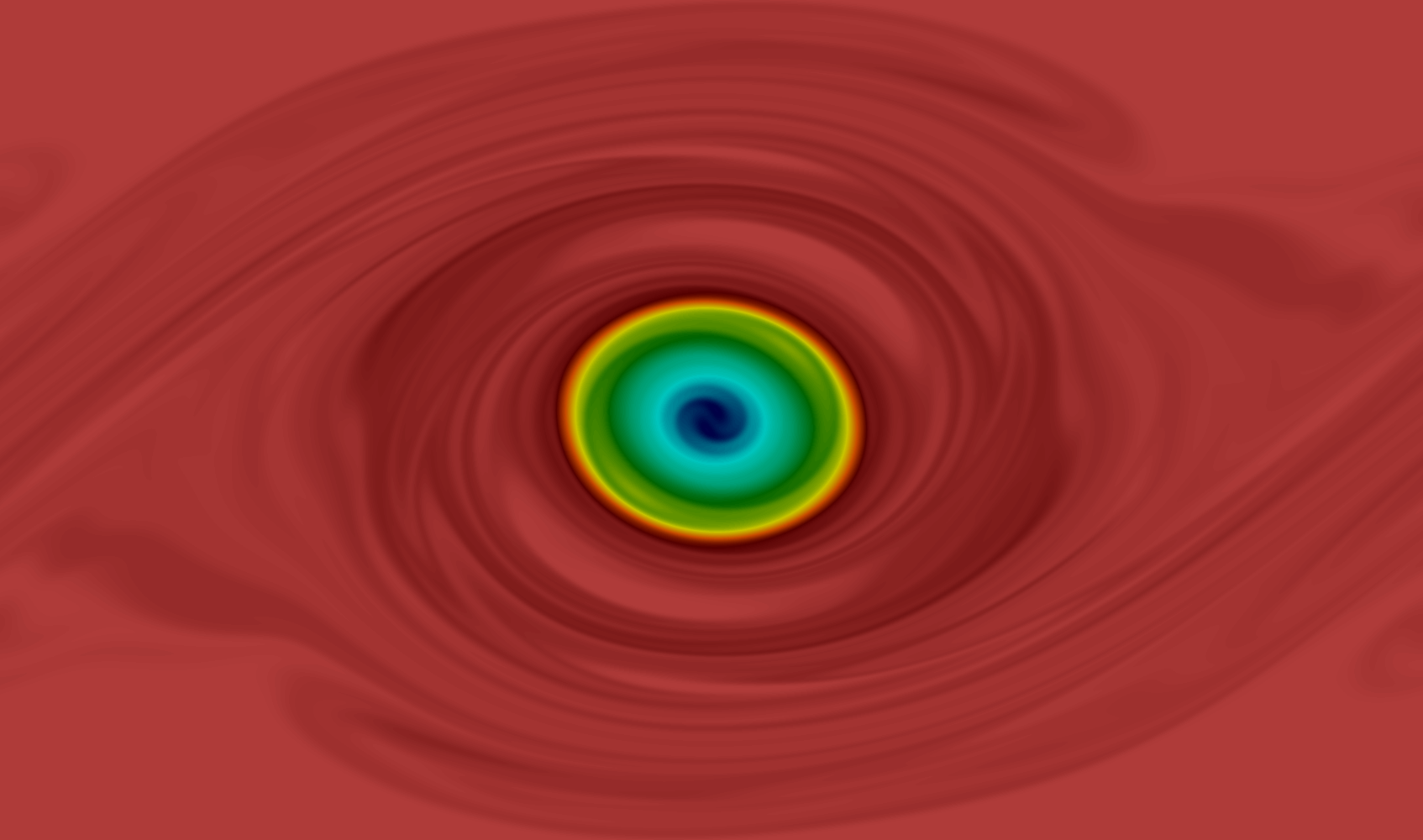}} 
\subfigure{\includegraphics[width=0.6\textwidth]{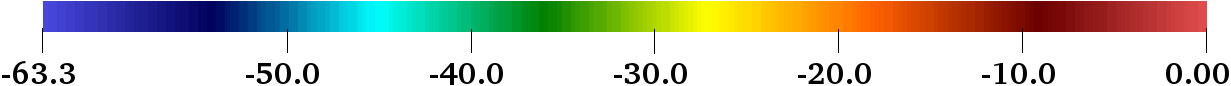}}
\caption{Vorticity $\nabla_h\times\uu_h\rb{t}$ for $\Rey=\num{10000}$ Kelvin--Helmholtz instability at (from left to right and top to bottom) $\bar{t}\in\set{5,10,17,34,56,200,240,278,400}$. Obtained with div-free $\HDIV$-FEM RT8 on the finest $256^2$ mesh; cf.\ Table \ref{tab:DOFs}.}
\label{fig:KH_Evolution_Vorticity}
\end{figure}

Let us begin with the description of our computational results in a rather qualitative way.
As is usual in the literature, plots of the vorticity are shown. 
Such a presentation, which illustrates the time evolution and dynamics of the involved vortices, can be seen in Figure~\ref{fig:KH_Evolution_Vorticity}.

%-----------------
%-----------------

Note that we show these results only for the highest Reynolds number $\Rey=\num{10000}$ and only on the finest $256^2$ mesh.
Thus, we only present results for our most resolved solution.
In the first row, the transition from the initial condition to the four primary vortices is shown.
At $\tbar=17$ the fine scales of the flow are clearly visible.
The four vertices are unstable in the sense that they have the tendency to merge.
This is a well-known property of two-dimensional flows for which (in contrast to 3D flows) energy is transferred from the small to the large scales.
In the second row, after the second merging process is completed at $\tbar=56$, we observe that the two vortices are rotating a very long time and are still clearly separated at $\tbar=200$. 
This is a very important difference compared to comparable computations in the literature.
To the best of our knowledge, until now there are no reliable results available in which the vortices are stable for such a long time.
To obtain these results, it was mandatory to resolve all the fine scales of the complex vortex dynamics. 
Also, one can see that directly after the merging process to two vortices, they have an ellipsoidal shape and fine scales are clearly visible.
As the two vortices rotate, shear forces act dissipatively on the fine scales, thereby smoothing and smearing them out which has the result that the shape of the vortices becomes more circular.
Finally, the two vortices start the pairing process in the third row.
Again, the resulting vortex at first has many fine scale details which get dissipated over time.
Note that the last vortex rotates in the middle of the domain.
However, as we will see later in more detail, predicting the time instance of the last vortex pairing cannot be done reliably due to the sensitivity of the problem.

\begin{figure}[t]
\centering
\subfigure{\includegraphics[width=0.32\textwidth]{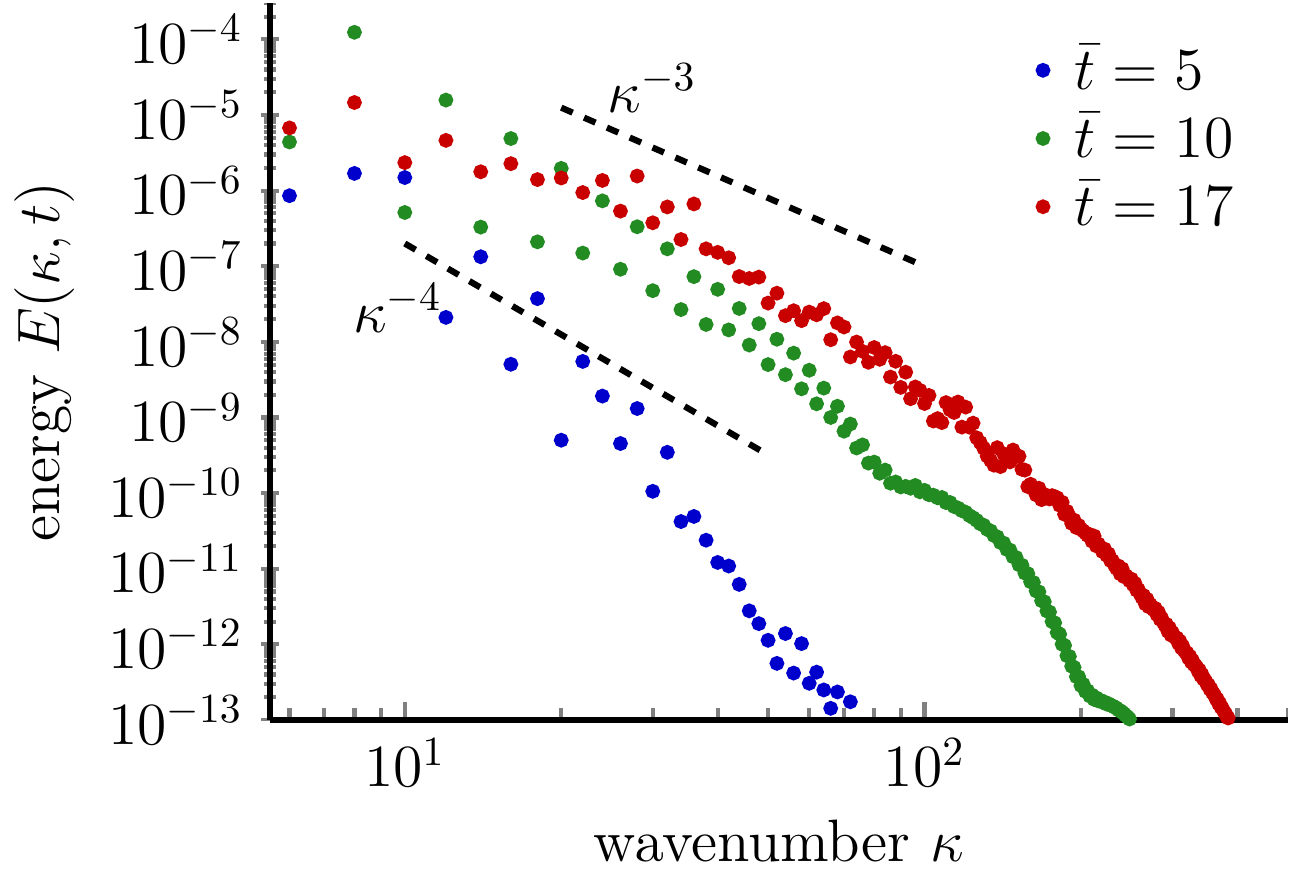}} \goodgap
\subfigure{\includegraphics[width=0.32\textwidth]{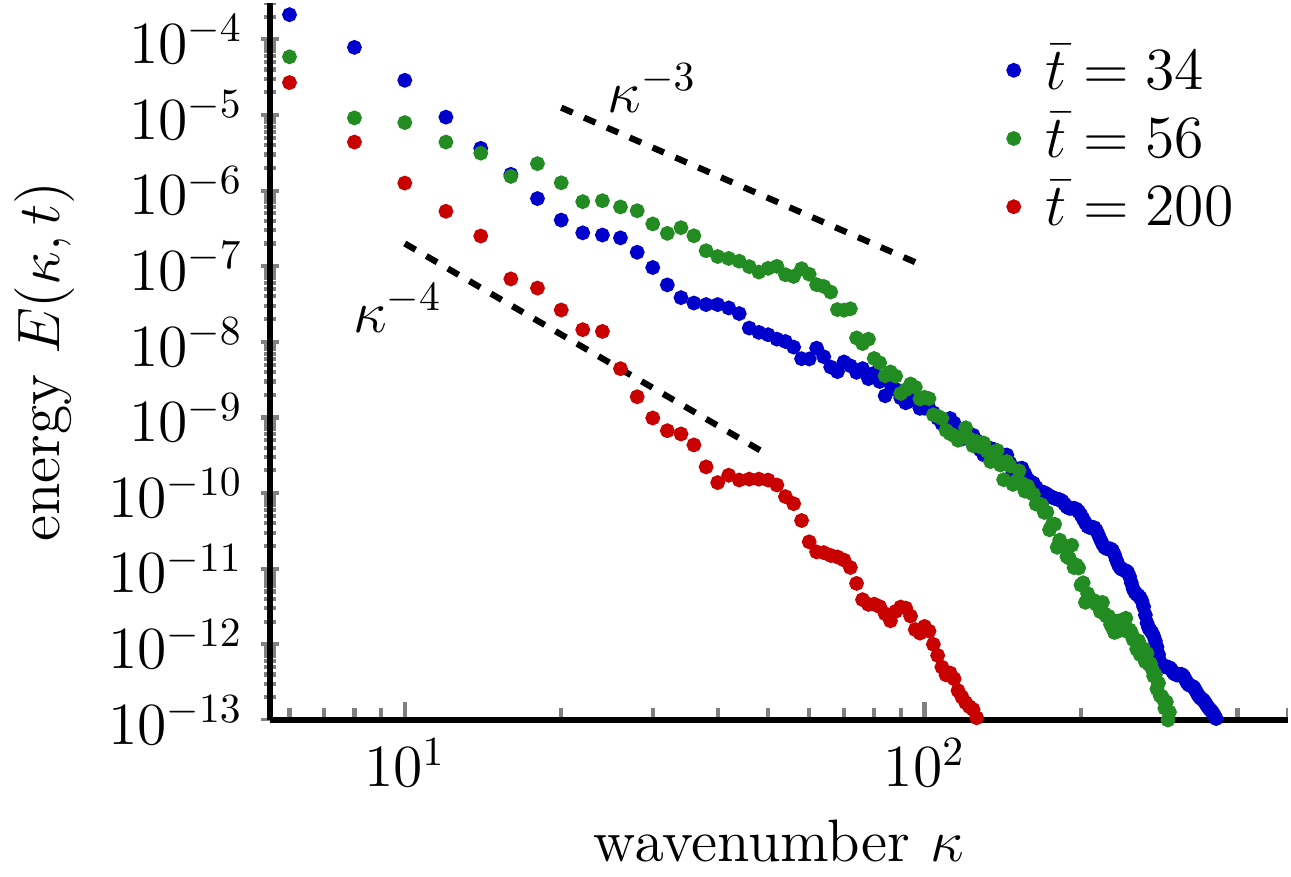}} \goodgap
\subfigure{\includegraphics[width=0.32\textwidth]{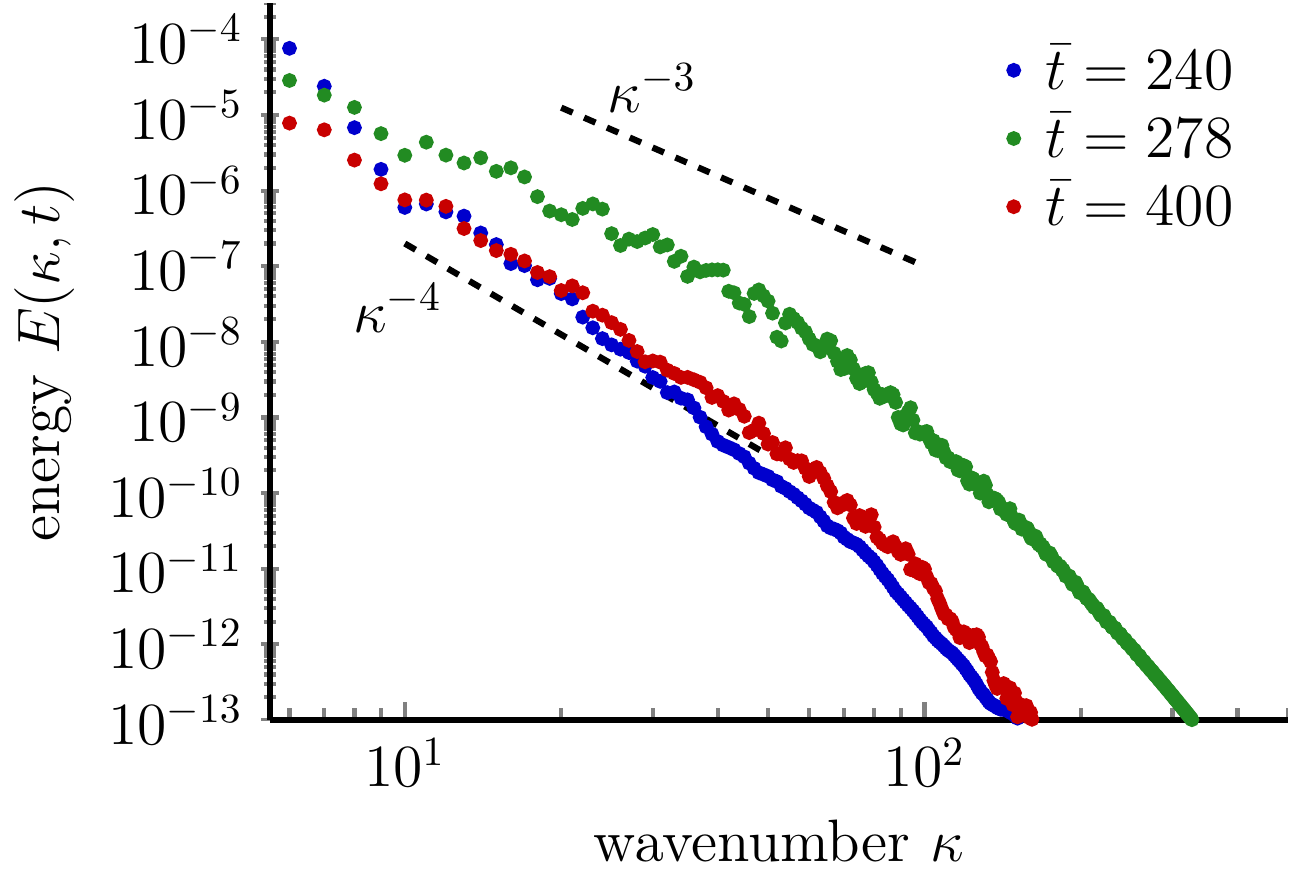}} 
\caption{Energy spectra $E\rb{\kappa,t}$ corresponding to Figure \ref{fig:KH_Evolution_Vorticity}. Each plot represents one row in Figure \ref{fig:KH_Evolution_Vorticity}.}
\label{fig:KH_Energy_Spectra}
\end{figure}

Accompanying these results, Figure~\ref{fig:KH_Energy_Spectra} shows the longitudinal energy spectra obtained at exactly the same points in time corresponding to Figure~\ref{fig:KH_Evolution_Vorticity}.
Each plot presents one row of the vorticity snapshots from Figure~\ref{fig:KH_Evolution_Vorticity}.
The first one shows that as the four primary vortices develop, more and more energy is inserted into the fine scales.
At $\tbar=17$, one observes the expected behavior for two-dimensional turbulence which says that the energy should be distributed according to something between $\kappa^{-3}$ and $\kappa^{-4}$.
Regarding the second plot, it can be observed that during the rotating period of two vortices, viscosity effects clearly remove energy from the fine scales (high wavenumbers/high frequencies); see the spectrum at $\tbar=200$.
However, immediately after the last pairing, the spectrum at $\tbar=278$, along with more energy in the fine scales, again shows the characteristics for 2D turbulence.

%-----------------
%-----------------

%---------------------------------------Kinetic energy-------------------------------------------------
\subsection{Kinetic energy}

%-----------------
\begin{figure}[t]
	\centering
	\includegraphics[width=0.85\textwidth]{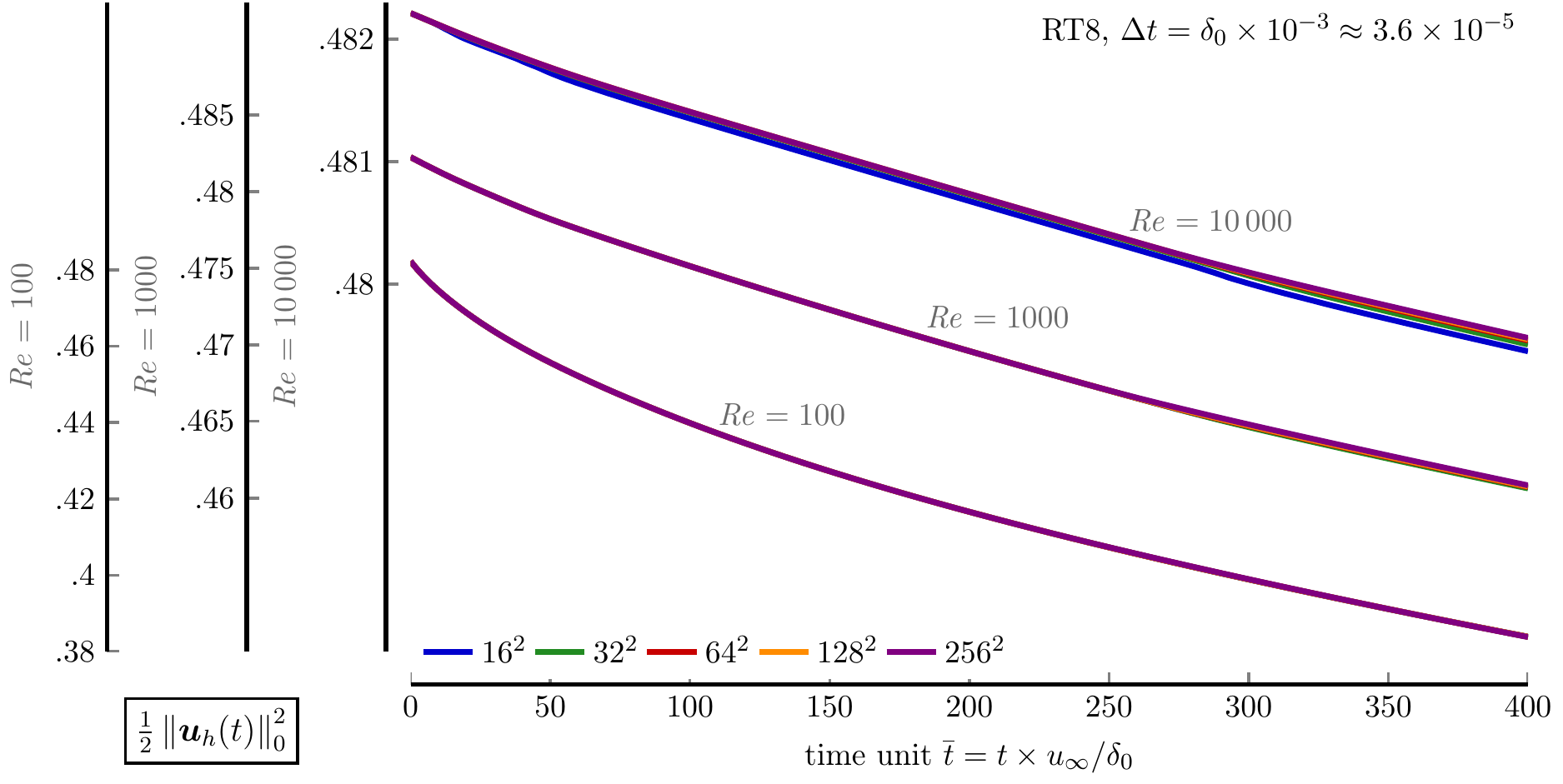}
	\caption{Evolution of kinetic energy $\Kin\rb{t,\uu_h}$ for different Reynolds numbers with divergence-free RT8 HDG and SBDF2 IMEX time stepping on a sequence of square meshes. Only for $\Rey=\num{10000}$, the mesh has a visible impact on the kinetic energy.}
	\label{fig:KineticEnergy}
\end{figure}
%-----------------

In Figure~\ref{fig:KineticEnergy}, the evolution of the resulting kinetic energy for all considered Reynolds numbers $\Rey\in\set{\num{100},\num{1000},\num{10000}}$ on all meshes $16^2,\dots,256^2$ can be seen.
Note that all families of curves start with the same amount of initial kinetic energy $\Kin\rb{0,\uu_h}\approx\num{0.4822}$.
For each Reynolds number, there is a different ordinate axis shown and thus, the offset in the curves is purely for a better view.

The first important observation is that the kinetic energy seems to be rather indifferent with respect to the mesh size.
Only for the highest Reynolds number, it is possible to detect a visible difference between the curves (note, however, the narrow scaling of the plot for $\Rey=\num{10000}$).
On the coarsest mesh ($16^2$), too much energy gets dissipated by artificial viscosity effects.
For $\Rey=\num{1000}$ and $\Rey=\num{100}$, the curves are basically congruent.
From our experience, kinetic energy is the quantity of interest which is easiest to compute accurately, even on coarse meshes. 

%-----------------
\begin{thmRem}
	In none of the existing literature concerning this problem \cite{GravemeierEtAl05,Burman07,AhmedEtAl17,YangBadiaCodina16} where it is seriously attempted to obtain a mesh-converged solution with respect to $\Kin$, the results are comparably convincing as in Figure~\ref{fig:KineticEnergy}.
	Only in \cite{SchroederLube18}, it is possible to observe mesh-convergence for the kinetic energy \textemdash{} however, the vortices merge much earlier there. 
\end{thmRem}
%-----------------

Secondly, we observe that regardless of the Reynolds number, all solutions have a strictly monotonically decreasing kinetic energy over time. 
Physically, this is the correct behavior as explained in Section \ref{sec:QuantitiesOfInterest}; see also equation \eqref{eq:EW}.
Due to a more dominant role of viscous effects, the flow loses comparably more kinetic energy for $\Rey=\num{100}$.
In fact, the results for $\Rey=\num{10000}$ show that the corresponding flow loses only about $\SI{0.54}{\percent}$ energy over the course of $\tbar=400$ time units compared to the energy at $t=0$.
For $\Rey=\num{100}$, on the other hand, the energy loss amounts to $\SI{20.41}{\percent}$ whereas for $\Rey=\num{1000}$, $\SI{4.42}{\percent}$ of the initial kinetic energy has been dissipated at $\tbar=400$.

%-----------------
\begin{thmRem}
	Note that for example in \cite{GravemeierEtAl05,AhmedEtAl17}, the kinetic energy seems to oscillate in time.
	This is clearly not a physical behavior and there is no mechanism in this problem which could excite such phenomenon.
	Therefore, the results presented in Figure~\ref{fig:KineticEnergy} are 
	physically more meaningful.
\end{thmRem}
%-----------------

%-------------------------Enstrophy and vorticity thickness----------------------------------------------
\subsection{Enstrophy and vorticity thickness}

%-----------------
\begin{figure}[t]
	\centering
	\includegraphics[width=0.85\textwidth]{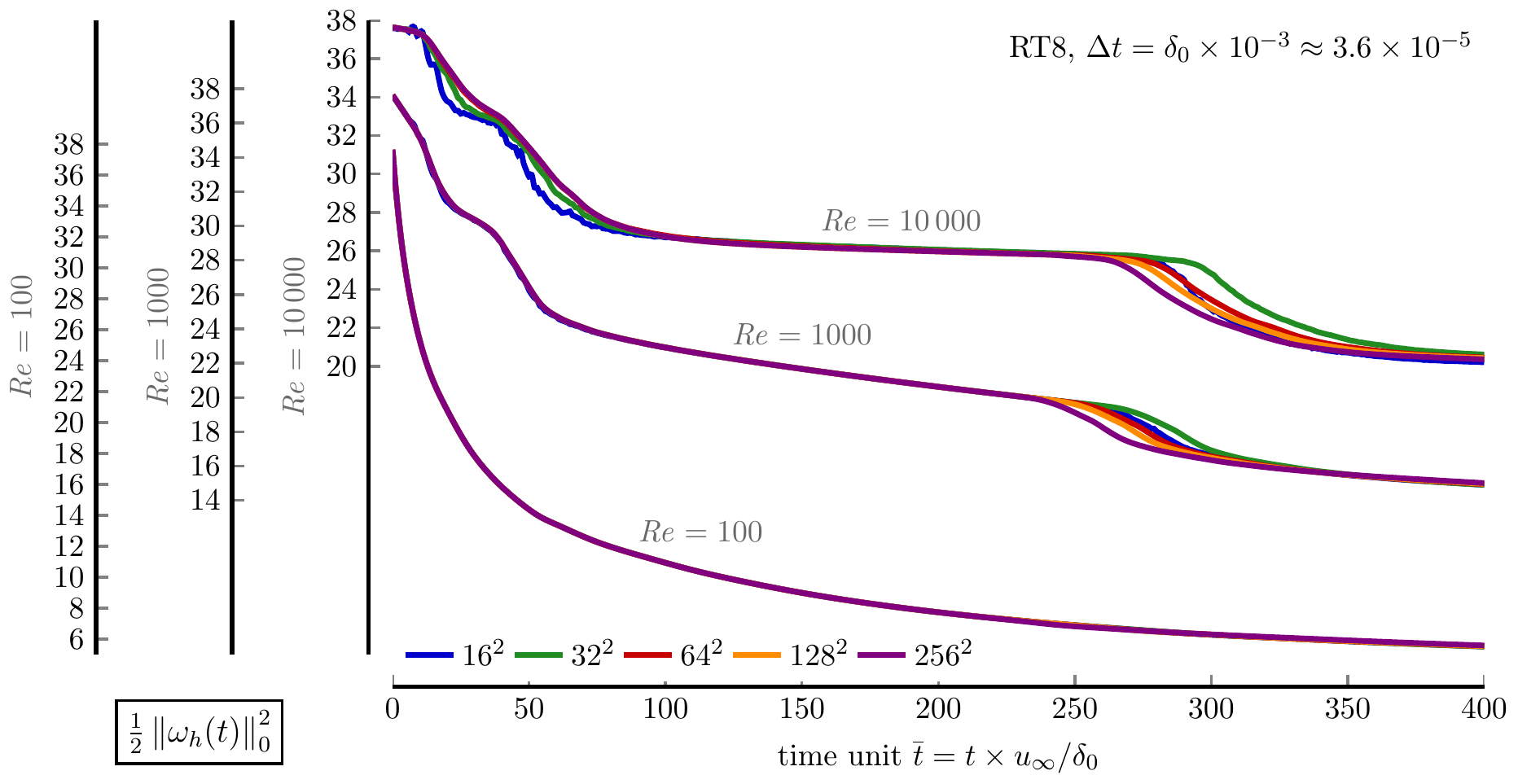}
	\caption{Evolution of enstrophy $\Ens\rb{t,\uu_h}$ for different Reynolds numbers with RT8 HDG and SBDF2 IMEX time stepping on a sequence of square meshes. Only for $\Rey=\num{100}$ the enstrophy is mesh-invariant at all times. }
	\label{fig:Enstrophy}
\end{figure}
%-----------------

Figure~\ref{fig:Enstrophy} presents the evolution of the  enstrophy for all considered Reynolds numbers on all meshes. 
At the initial time, one has roughly the same amount of enstrophy $\Ens\rb{0,\uu_h}\approx\num{37.63}$ for all simulations.
Again, for each Reynolds number, there is a different ordinate axis shown and thus, the offset in the curves is purely for a better view.

In contrast to the kinetic energy, different meshes result in a different enstrophy already for $\Rey=\num{1000}$. 
However, the only real difference can be seen in the interval $\tbar\in\sqb{250,325}$ where the last pairing of the vortices takes place.
For $\Rey=\num{10000}$, the enstrophy behaves differently on the coarse meshes for the first two pairings at the beginning of the simulation.
Nevertheless, Figure~\ref{fig:Enstrophy} shows mesh convergence at least for $\tbar<200$ for the three finest meshes.
Only for $\Rey=\num{100}$, all results are mesh-independent at all times.
At least up to $\tbar=200$, we believe that our results are reliable for all Reynolds numbers.
For $\Rey\in\set{\num{1000},\num{10000}}$, the last merging from two vortices to one vortex seems to be very sensitive with respect to perturbations to such an extent that even the highest solution fails to show mesh-independence.

Qualitatively, with the possible exception of the very coarse $16^2$ simulation for $\Rey=\num{10000}$, all curves are strictly monotonically decreasing.
This is in agreement with the theoretical considerations in Section~\ref{sec:SelfOrga2DTurbulence}.
Indeed, equation \eqref{eq:EW} predicts that the decrease in enstrophy is especially strong whenever the palinstrophy is high and, in anticipation of Figure \ref{fig:Palinstrophy}, this behavior is very well observable for $\Rey\in\set{\num{1000},\num{10000}}$.

%-----------------
\begin{thmRem}
Note that the desire of having a mesh-converged enstrophy \emph{pointwise in time} is rather ambitious and, in the present situation, not necessarily realistic from an analytical point of view.
However, exclusively in the 2D periodic case, one can at least find $\Lp{\infty}{\rb{\HM{1}{}}}$ \emph{energy} estimates for the exact solution $\uu$; see, for example, \cite[Section 5.4]{DoeringGibbon95}.
Such energy estimates can be extended to FE approximations $\uu_h$ but the question of having control over $\norm{\uu-\uu_h}_\Lp{\infty}{\rb{\HM{1}{}}}$ by means of \emph{error} estimates with quantitative convergence rates is beyond the scope of this paper and the authors are not aware of any existing literature in this direction.
In fact, numerical error analysis for time-dependent Navier--Stokes flows usually only covers the convergence of the integral quantity $\int_0^t \norm{\nabla_h\sqb{\uu-\uu_h}\rb{\tau}}_0^2\dtau$ for $t\in\rb{0,\tend}$; cf., for example, \cite{John16,SchroederEtAl18}.
By choosing a divergence-free method, also on the discrete level, the $\LP{2}{}$-norm of the gradient is identical to the $\LP{2}{}$-norm of the vorticity \textemdash{} the enstrophy.
Thus, divergence-free methods allow one to directly consider the error in the enstrophy which is, compare Section~\ref{sec:SelfOrga2DTurbulence}, a very important quantity in 2D turbulence.
If the results from Figure~\ref{fig:Enstrophy} are post-processed to show the \emph{time integral} of the enstrophy, convergence can indeed be observed for our results (for brevity not shown here).
\end{thmRem}
%-----------------

%-----------------
\begin{figure}[t]
	\centering
	\includegraphics[width=0.85\textwidth]{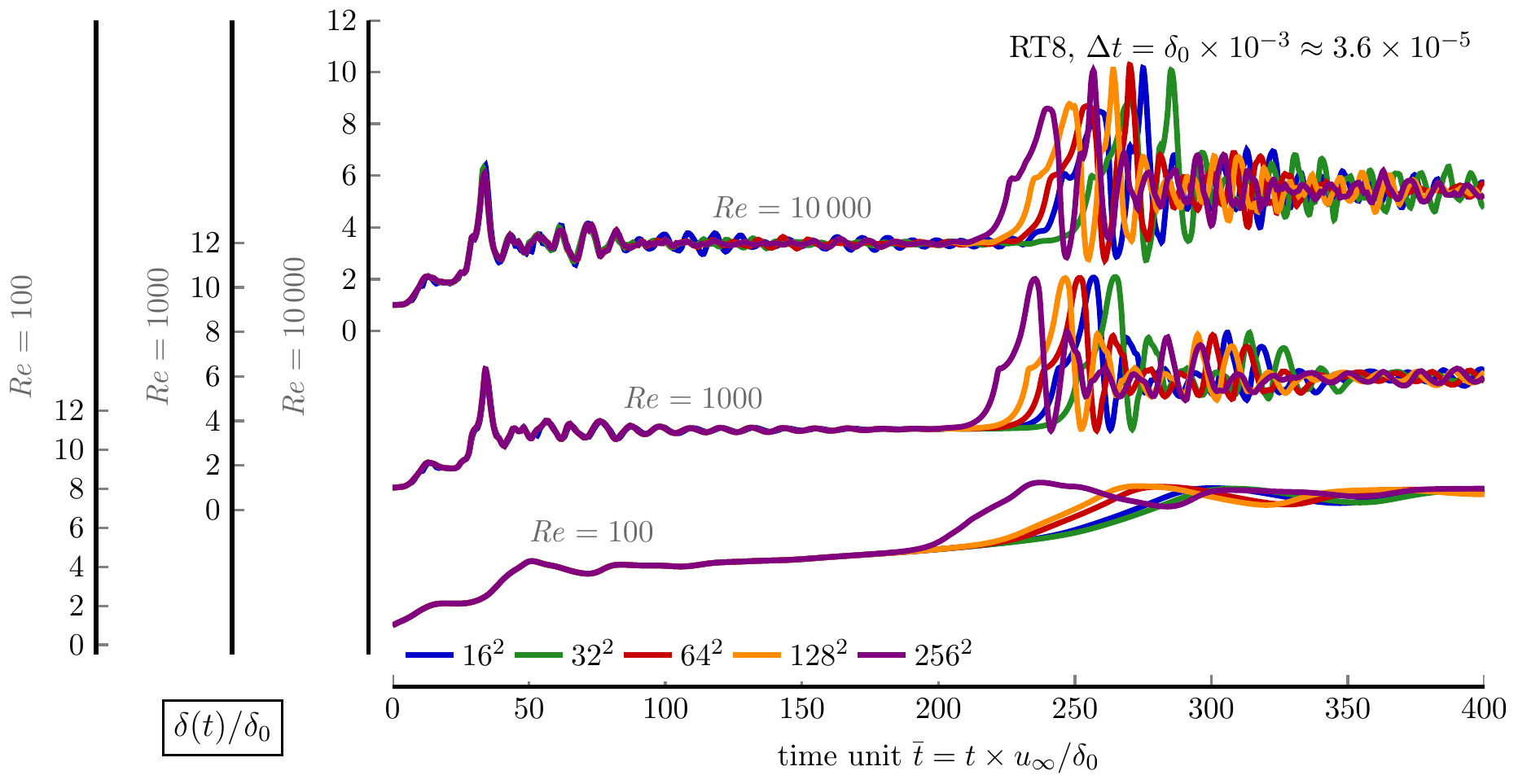}
	\caption{Evolution of scaled vorticity thickness $\delta\rb{t}/\delta_0$ for different Reynolds numbers with RT8 HDG and SBDF2 IMEX time stepping on a sequence of square meshes. }
	\label{fig:VorticityThickness}
\end{figure}
%-----------------

In Figure~\ref{fig:VorticityThickness}, the evolution of the resulting relative vorticity thickness for all considered Reynolds numbers $\Rey\in\set{\num{100},\num{1000},\num{10000}}$ on all meshes $16^2,\dots,256^2$ can be seen.
For $\Rey\in\set{\num{1000},\num{10000}}$, the local maxima at $\tbar\approx 34$ correspond to the pairing process from four to two vortices.
This local maximum is less pronounced for $\Rey=\num{100}$ and occurs later at $\tbar\approx 51$.
Afterwards, corresponding to an ellipsoidal shape of the vortices, $\delta\rb{t}/\delta_0$ oscillates until the next pairing is imminent.
We note that for all Reynolds numbers, the plots show mesh-convergence up to $\approx\num{200}{\tbar}$ time units, which means that the pairing processes can be predicted reliably on all meshes right to the point where the last two vortices merge.
A mesh-independent capturing of the last merging, which is most probably very sensitive to perturbations, could not be achieved. 
Then, the last merging takes place after which the vorticity thickness again oscillates.
As time proceeds, the last resulting vortex becomes more and more circular.

%---------------------------------------Palinstrophy-------------------------------------------------
\subsection{Palinstrophy}

\begin{figure}[t]
	\centering
	\includegraphics[width=0.85\textwidth]{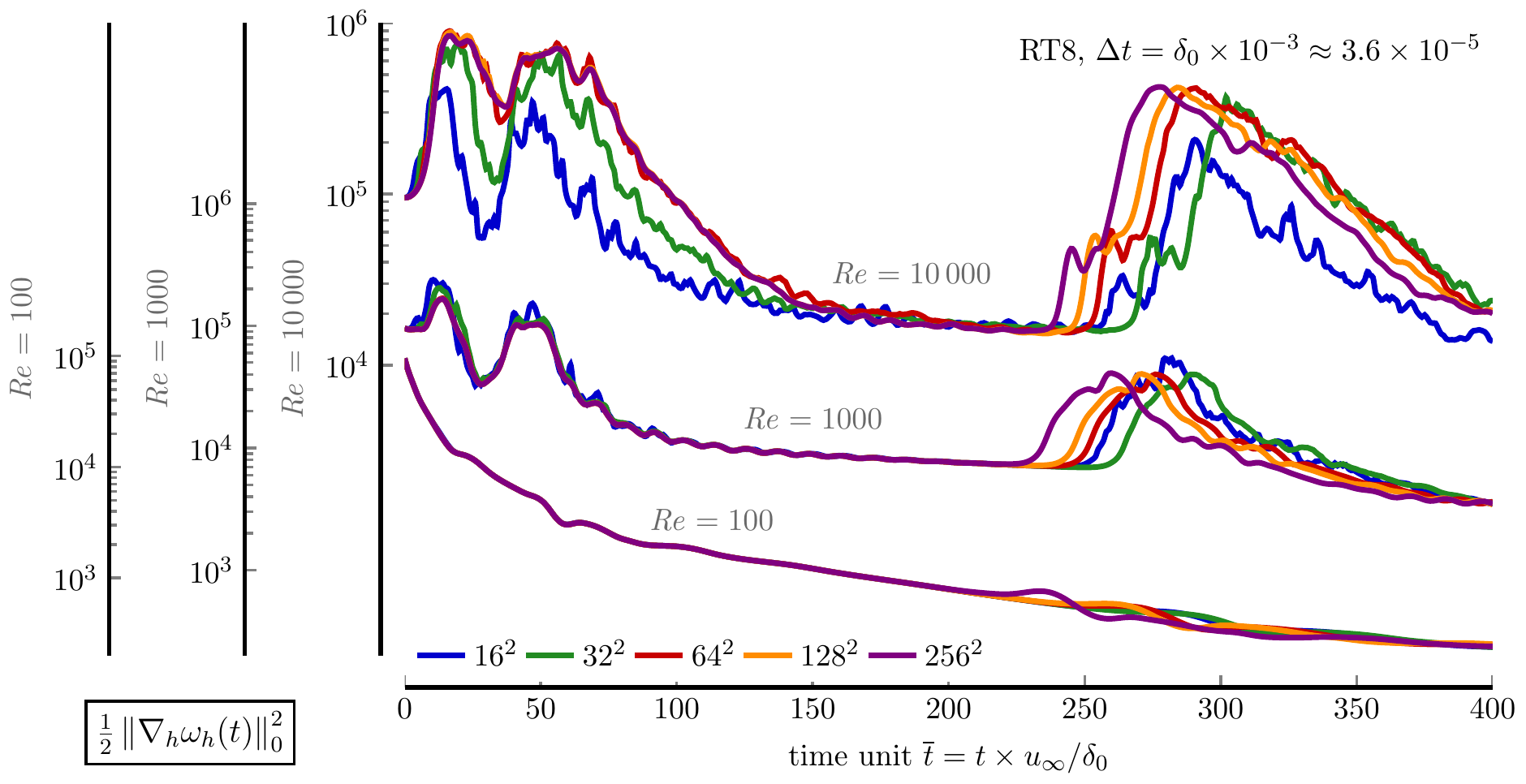}
	\caption{Evolution of palinstrophy $\Pal\rb{t,\uu_h}$ for different Reynolds numbers with RT8 HDG and SBDF2 IMEX time stepping on a sequence of square meshes. The palinstrophy is not even mesh-converged at all times for the smallest Reynolds number.}
	\label{fig:Palinstrophy}
\end{figure}

The evolution of the palinstrophy for all Reynolds numbers and on all meshes is 
displayed in Figure~\ref{fig:Palinstrophy}. The initial palinstrophy is given 
by $\Pal\rb{0,\uu_h}\approx\num{95219}$.
Once more, for each Reynolds number, there is a different ordinate axis shown and thus, the offset in the curves is just for a better view.

One can observe that for $\Rey=\num{100}$ the palinstrophy starts to be mesh-dependent beginning from $\tbar=200$. 
For $\Rey=\num{1000}$ and, more dramatically, for $\Rey=\num{10000}$ the palinstrophy is obviously a very sensitive quantity of interest which makes it perfect for comparing results.
In contrast to kinetic energy and enstrophy, palinstrophy can increase spontaneously.
Such outbursts always correspond to the merging of vortices in the Kelvin--Helmholtz problem.
Especially for $\Rey=\num{10000}$, one can see three time intervals where the palinstrophy has very pronounced local maxima. 
These intervals mark the three merging processes of the vortices.
Especially, this means that independent of the Reynolds number, the last merging processes does not occur before $\tbar=200$.
On other hand, the time instance of the last pairing cannot be predicted precisely which is due to the sensitivity of the problem.
Note, however, that although the magnitude of the palinstrophy is strongly mesh-dependent, the points in time where the first two pairings occur can be approximately identified independently of the particular mesh resolution.
Consequently, also for all Reynolds numbers, our results are mesh-independent for $\tbar<200$ which makes them reliable there.

%-----------------
%-----------------

%-----------------
\begin{thmRem}
One can evidently see that it is very hard to reliably compute approximate solutions for this problem for times larger than $\tbar=200$.
In fact, we want to remark at this point	 that also numerical error analysis predicts an \emph{exponential} growth in the error for time-dependent Navier--Stokes simulations; cf., for example \cite{John16}.
More precisely we refer to \cite{SchroederEtAl18} where this point is both proven and discussed for the non-hybridized variant of the space discretization introduced in Section~\ref{sec:DivFreeHDG}.
There, it is shown that, roughly speaking, the error in kinetic energy $\Kin$ and enstrophy $\Ens$ can only be controlled up to a factor $\exp\rb{G_\uu\rb{t}}$, where $G_\uu$ is a Gronwall term which depends on the regularity of the exact solution $\uu$.
As $t\to\infty$, unavoidably, one loses control over the accuracy of any finite element approximation.
There is also the numerical example of the `2D planar lattice flow', which consists of four counter-rotating vortices and shares some structural properties of the Kelvin--Helmholtz problem.
In a simulation of the planar lattice flow, it turned out that any numerical method can uphold the structure of those vortices for a certain time \textemdash{} however, eventually, the structures collapse.
This behavior might also facilitate the understanding of the difficulties in obtaining a reliable last merging for the Kelvin--Helmholtz problem. 
\end{thmRem}
%-----------------

\begin{figure}[t]
	\centering
	\includegraphics[width=0.49\textwidth]{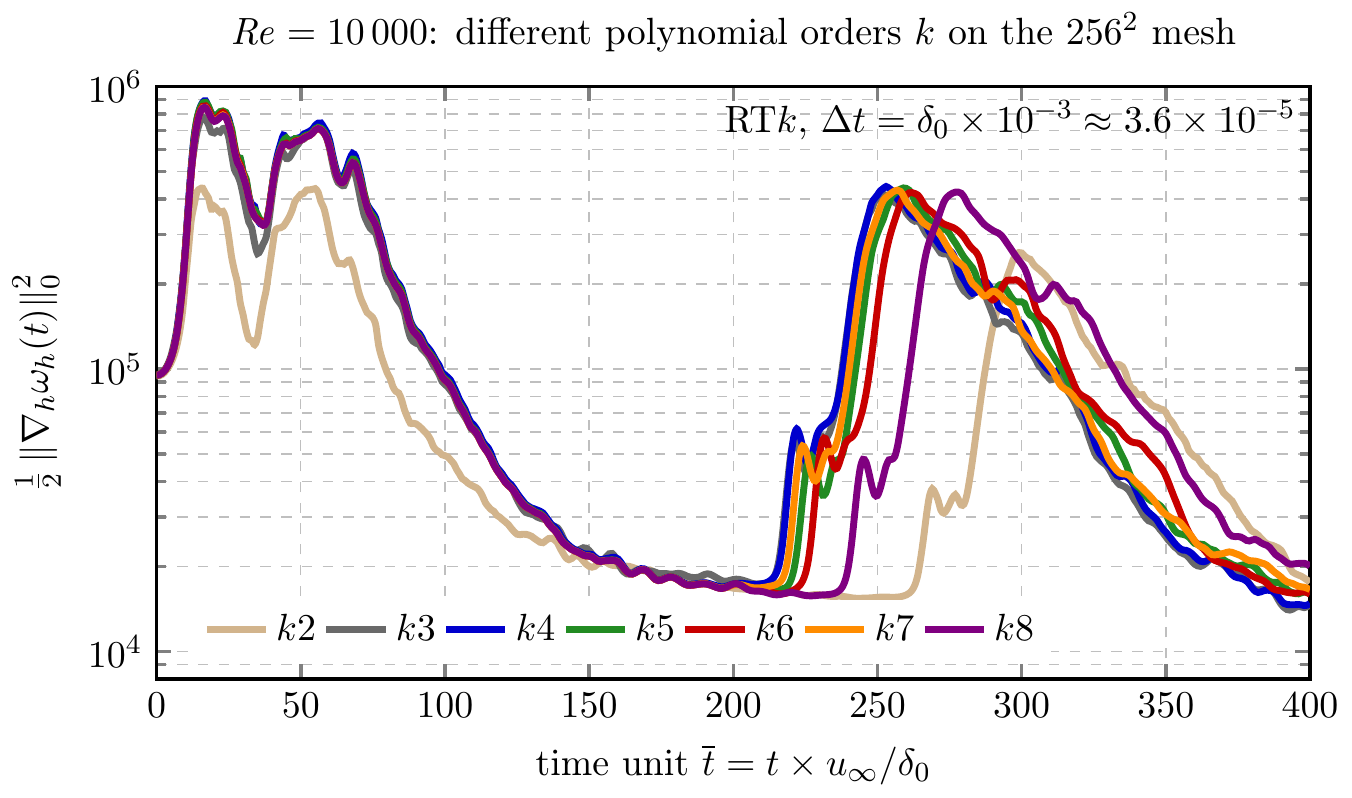} 
	\caption{Evolution of palinstrophy $\Pal\rb{t,\uu_h}$ for $\Rey=\num{10000}$ with RT$k$ HDG and SBDF2 IMEX time stepping on the $256^2$ mesh with varying polynomial degree $k\in\set{2,3,4,5,6,7,8}$.}
	\label{fig:Palinstrophy_kConv}
\end{figure}

Furthermore, in order to also consider the behavior of the Kelvin--Helmholtz instability problem in the situation that the mesh is fixed and the polynomial order $k$ of the discrete spaces is changed, we  provide Figure~\ref{fig:Palinstrophy_kConv}.
This graph shows the evolution of the palinstrophy, which previously has been identified as the most sensitive quantity of interest, as a function of $k$ where the finest $256^2$ mesh has been chosen.
Firstly, $k$-convergence of $\Pal\rb{t,\uu_h}$ can be observed for $\tbar\leqslant 200$ which is in agreement with the $h$-convergence study above.
Furthermore, $k=4$ on the $256^2$ seems to be the minimum resolution for which reliable results can be computed up to $\tbar=200$.
For larger times $\tbar > 200$ the palinstrophy behaves analogously as for the $h$-convergence study.
Namely, the instance in time where the last merging occurs is also very sensitive with respect to the used polynomial degree of the FE spaces.

%-----------------
%-----------------

%--------------------------------- Numerical dissipation ---------------------------------------------
\subsection{Numerical dissipation}

In the previous section, we investigated the convergence behavior of kinetic energy, enstrophy, vorticity thickness and palinstrophy under mesh (and partially $k$) refinement.
It turned out that our results are reliable for $\tbar\leqslant 200$. 
In order to confirm this statement, we additionally want to investigate the numerical dissipation induced by the discretization scheme.
Note that numerical dissipation is a quantity of interest which is relevant for evaluating the accuracy of the numerical method itself.
In this sense, and opposed to the quantities of interest presented in Section~\ref{sec:QuantitiesOfInterest}, numerical dissipation does not characterize the flow problem at hand. 
This is why this subsection is clearly separated from the discussion of the characteristics of the Kelvin--Helmholtz instability problem.
However, numerical dissipation is frequently investigated and satisfactory results for a DNS flow simulation can only be obtained when it is comparatively small.

The basic idea is to investigate how \eqref{eq:EW} behaves in the discrete setting and thus, the time-dependent numerical dissipation $\eps_h\rb{t}$ is frequently defined as \cite{GassnerBeck13,FehnWallKronbichler2018} 
\begin{equation}
	\frac{\drm}{\drm t}\Kin{\rb{t, \uu_h}}
		= -2\nu\Ens{\rb{t,\uu_h}} - \eps_h\rb{t}. 
\end{equation} 
Then, from this equality, by integration over $\sqb{0,t}$ for all $t\in\rsb{0,\tend}$, an absolute and a relative (w.r.t.\ the loss of kinetic energy) integrated numerical dissipation can be defined as
\begin{equation}
	\eps_h^\intrm\rb{t}
		= \abs{\Kin\rb{0, \uu_h}-\Kin\rb{t, \uu_h}-2\nu\int_0^t \Ens\rb{\tau, \uu_h} \dtau}, \quad
		\eps_{h,\rel}^\intrm\rb{t}
			= \frac{\eps_h^\intrm}{\abs{\Kin\rb{0, \uu_h}-\Kin\rb{t, \uu_h}}}, 
\end{equation}
which measures the total amount of numerical dissipation present in $\sqb{0,t}$.
In Table~\ref{tab:NumDiss}, both quantities are shown for the most difficult $\Rey=\num{10000}$ computations, where the integration has been performed in $\tbar\in\sqb{0,200}$. 
One can see that the relative numerical dissipation is very small, whereas the absolute numerical dissipation is already close to zero on the finest mesh.

%-----------------
\begin{table}[t]
\caption{Numerical dissipation (absolute and relative) with RT8 HDG and SBDF2 IMEX time stepping for $\Rey=\num{10000}$ on the considered sequence of square meshes.}
\label{tab:NumDiss}
\centering 
\begin{tabular}{crrrrr} 
\toprule
	Mesh	
		& $16^2$ 
		& $32^2$		
		& $64^2$	
		& $128^2$
		& $256^2$
		\\ 
\otoprule
	$\eps_h^\intrm\rb{\tbar=200}$		
		& \pgfmathprintnumber[precision=1]{7.608516363643e-05}	
		& \pgfmathprintnumber[precision=1]{1.452210862242e-05}	
		& \pgfmathprintnumber[precision=1]{1.422335464224e-06}	
		& \pgfmathprintnumber[precision=1]{3.417780264619e-08}	
		& \pgfmathprintnumber[precision=1]{3.487552908575e-09}	
		\\
	$\eps_{h,\rel}^\intrm\rb{\tbar=200}$		
		& \pgfmathprintnumber[precision=1]{4.959817026347e-02}		
		& \pgfmathprintnumber[precision=1]{9.786371971634e-03}	
		& \pgfmathprintnumber[precision=1]{9.631268419139e-04}	
		& \pgfmathprintnumber[precision=1]{2.317081358294e-05}	 
		& \pgfmathprintnumber[precision=1]{2.364319006220e-06}	
		\\
\bottomrule		
\end{tabular}
\end{table}
%-----------------

%------------------------------------------------------------------------------------------------
%------------------------------------SOURCES PERTURBATION----------------------------------------
%------------------------------------------------------------------------------------------------
\section{Computational studies of some possible sources for perturbations}
\label{sec:Perturbations}

In this section, we want to emphasize how sensitive the Kelvin--Helmholtz instability problem is with respect to some kinds of perturbation a numerical method can induce.
Thus, we consider as examples some possible sources of error.
Note that the sensitivity of this two-dimensional flow problem is strongly related to our explanations in Section~\ref{sec:SelfOrga2DTurbulence}.
In addition, this section is intended as a possible explanation why the results obtained for this problem vary so extremely in the existing literature \textemdash{} Kelvin--Helmholtz instabilities are very fragile flow systems and therefore hard to compute reliably.
We only show the evolution of the palinstrophy in this section as it has clearly emerged as the  most sensitive quantity of interest.

Furthermore, as seen in the Van Groesen theory for the continuous problem in Section~\ref{sec:SelfOrga2DTurbulence}, potential perturbations $\xi \in E_1 \cup E_2 \cup \ldots \cup E_{k-1}$ cause a jump to a coarser invariant set $I_j$ with $j<k$ after some time. 
In numerical approximations, such perturbations may result from different sources like coarse meshes, insufficient numerical integration or resolution of the arising linear systems. 
Regarding the statement of (arbitrarily small) perturbations from sufficiently coarse 
eigenspaces (for the continuous problem), one has to be aware that the numerical error accumulates in time due to the deterministically chaotic character of the Navier--Stokes solution; cf. \cite{SchroederEtAl18}.
We emphasize that this section represents a sensitivity analysis aimed at investigating how perturbations impact the last pairing process from two vortices to one vortex.

%--------------------------------------Triangular meshes---------------------------------------
\subsection{Structured and unstructured triangular meshes}

Whenever a finite element method is used, maybe the first design decision is whether triangular or rectangular meshes (in 2D) are applied.
While this is, at least in our opinion, mostly a matter of taste, we show that for the Kelvin--Helmholtz instability problem, it can actually lead to completely different merging times of the involved vortices.
Furthermore, it makes a huge difference if structured or unstructured triangular meshes are used for the computation.

%-----------------
\begin{figure}[t]
	\centering
	\includegraphics[width=0.49\textwidth]{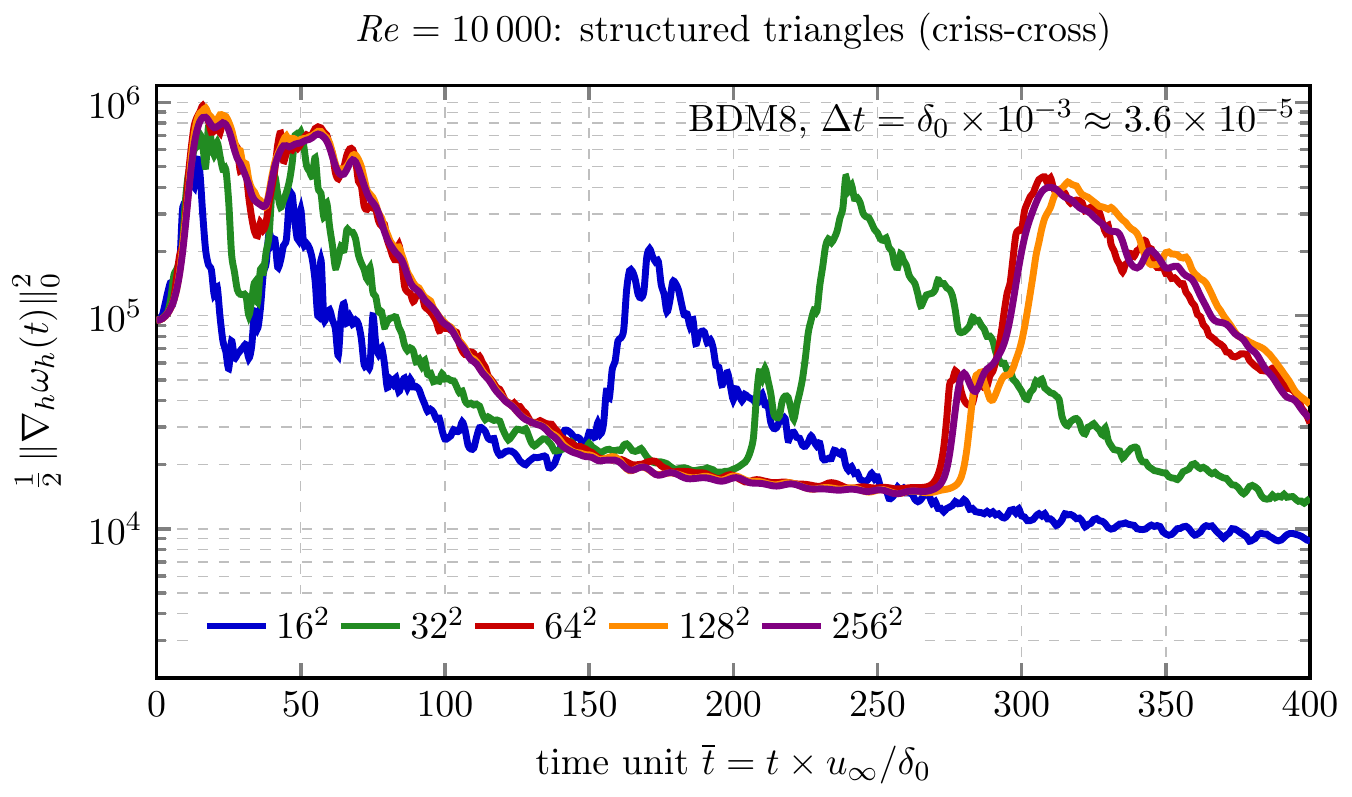} \goodgap
	\includegraphics[width=0.49\textwidth]{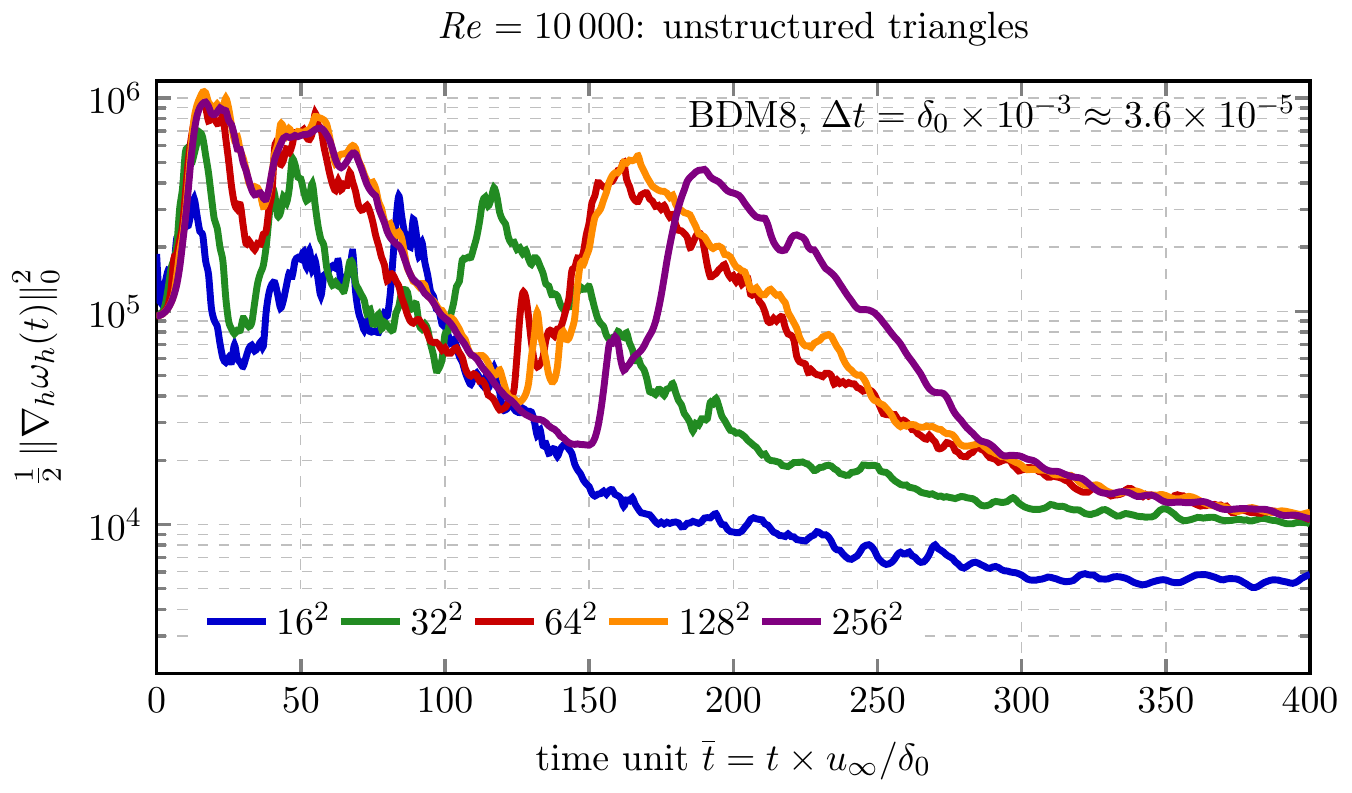}
	\caption{Evolution of palinstrophy $\Pal\rb{t,\uu_h}$ for $\Rey=\num{10000}$ with BDM8 HDG and SBDF2 IMEX time stepping on a sequence of structured criss-cross (left) and unstructured (right) triangular meshes.}
	\label{fig:Palinstrophy_TRIGS_BDM8}
\end{figure}
%-----------------

Thus, here, we repeat the computations of Section \ref{sec:ComputationalStudies} but change from meshes consisting of squares to meshes consisting of either structured or unstructured triangles.
In doing so, we continue to refer to the triangular meshes as $16^2$, $32^2$, etc.\ but in fact, this only means that the triangular meshes are chosen such that, very roughly, a comparable number of DOFs as for the corresponding square mesh (see Table~\ref{tab:DOFs}) are used.
More precisely, the structured triangular mesh results from the square mesh by dividing each patch of four squares into eight triangles, via the diagonals of the quadratic macro element.
  The unstructured triangular mesh, on the other hand, is generated by starting with a coarse triangular decomposition obtained by the mesh generator \texttt{netgen} \cite{schoberl1997netgen} (advancing front) which is then refined using edge bisection such that a comparable number of DOFs is obtained.
Note, however, that this section is not intended to give any recommendations concerning the choice of meshes or elements.
We only want to demonstrate the impact of this choice on the computational results. 

In Figure~\ref{fig:Palinstrophy_TRIGS_BDM8}, the evolution of the palinstrophy for the problem with $\Rey=\num{10000}$ is shown on two sequences of structured and unstructured triangular meshes, respectively.
One can see that as the structured mesh is refined (left), the last pairing of the two remaining vortices is delayed further and further.
In fact, on the finest mesh, the last merging happens even later as on the square meshes in Section~\ref{sec:ComputationalStudies}.

For unstructured meshes (Figure \ref{fig:Palinstrophy_TRIGS_BDM8}, right), the situation is slightly different.
On the coarsest mesh, the last merging happens extremely early, which is a situation that can often be found in the literature.
As the mesh is refined, this last merging process again shifts towards the end of the time interval.

%------------------------------ Linear systems & numerical integration ----------------------------------
\subsection{Inaccurate solution of linear systems and numerical integration} 
\label{sec:InaccurateSolution}

Let us now turn to more subtle, but nonetheless very common possible sources for perturbation.
In Section~\ref{sec:DivFreeHDG}, it became clear that we have to solve a large
linear system of the form $M^\ast \mathbf{x}=\mathbf{b}$ at every time step.
As explained in more detail in Section~\ref{sec:TimeDiscretization}, in order to compute our solutions in Section~\ref{sec:ComputationalStudies}, we used a sparse direct solver together with iterative refinement which ensures that linear systems are solved accurately up to a (relative) tolerance of $\num{e-12}$ (measured in the 2-norm).

\begin{figure}[t]
	\centering
	\includegraphics[width=0.49\textwidth]{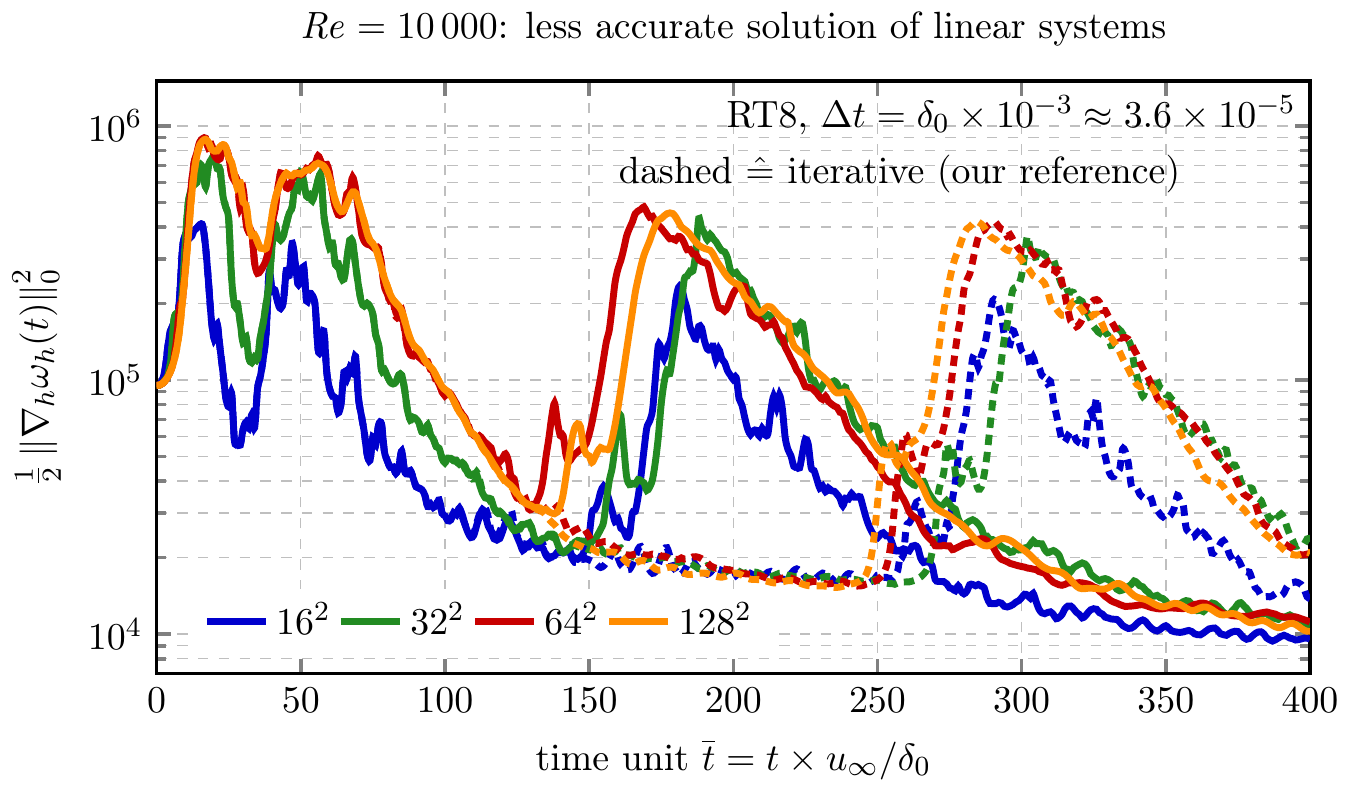} \goodgap
	\includegraphics[width=0.49\textwidth]{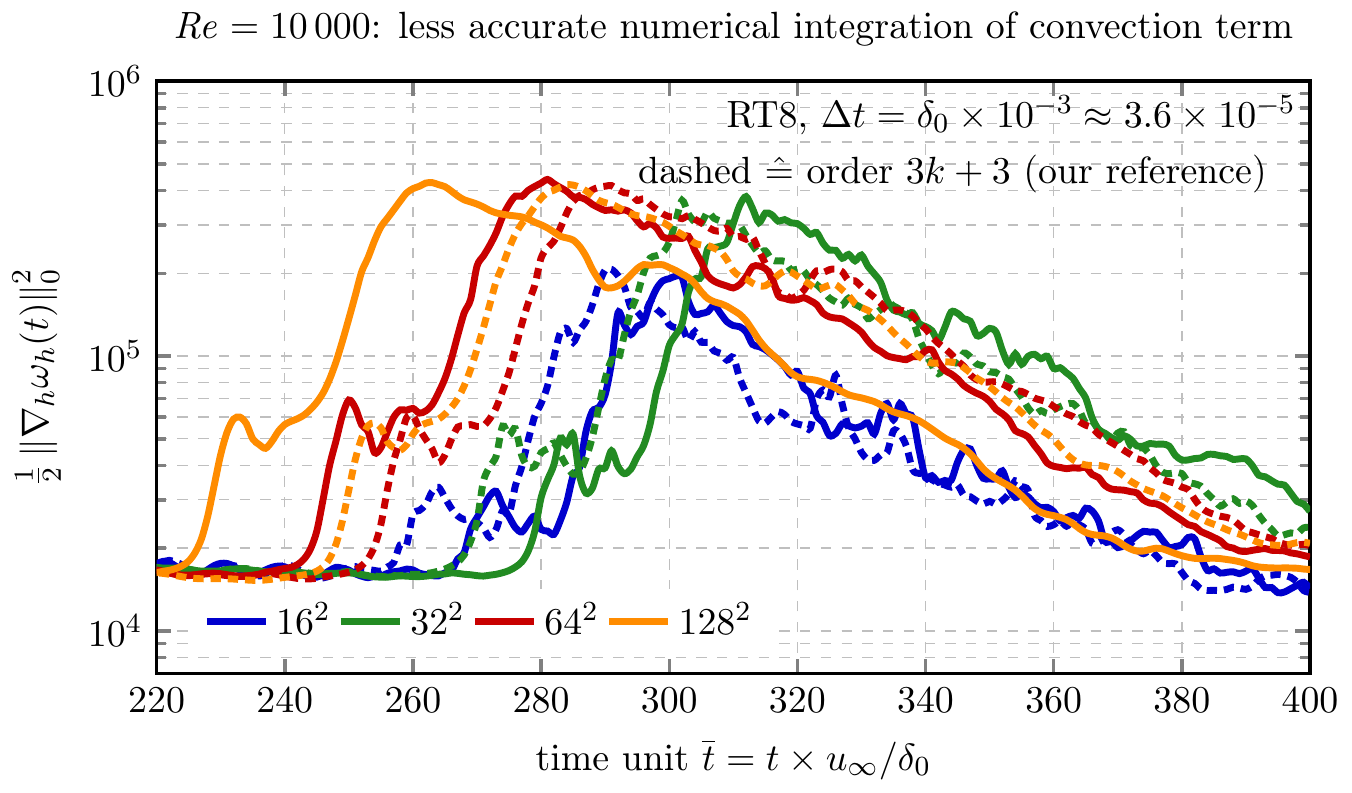}
	\caption{Evolution of palinstrophy $\Pal\rb{t,\uu_h}$ for $\Rey=\num{10000}$ with RT8 HDG and SBDF2 IMEX time stepping on a sequence of square meshes. Perturbation due to inaccurate solution of linear systems (left) and under-integration (right). Dashed line shows the corresponding solutions from Section \ref{sec:ComputationalStudies}. Right: For $\tbar < 220$, there is no visible difference. }
	\label{fig:Palinstrophy_INACCURATE_RT8}
\end{figure}

An alternative approach, probably used more frequently, consists in applying only the sparse direct solver (without iterative refinement) for solving $M^\ast \mathbf{x}= \mathbf{b}$. 
Here, we add a small pressure mass matrix scaled with $\num{-e-12}$ and apply the sparse Cholesky solver from \texttt{NGSolve} \cite{Schoeberl14}. 
Figure~\ref{fig:Palinstrophy_INACCURATE_RT8} (left) compares results from Section~\ref{sec:ComputationalStudies} with results that were obtained if the iterative scheme is replaced by the direct approach, keeping all other settings fixed.
We notice that applying the direct solver introduces a small perturbation due to the perturbation in the matrix, round-off errors, and condition numbers considerably larger than one. 
These small perturbations are significantly larger than the $\num{e-12}$ error tolerance guaranteed by the iterative solution. 
As a result, one observes that the merging of the last two vortices happens approximately $100\tbar$ time units earlier.
We believe that this sensitivity to accurate solutions of linear systems could also be a reason for partially completely different merging times in the literature.

%-----------------
%-----------------

A second source of error when using numerical schemes is the issue of numerical integration; more precisely, let us briefly investigate the effects of under-integration; cf., for example, \cite{BeckEtAl16,Kopriva17}.
For example, we consider the accuracy of the numerical integration of the convection term, see Section~\ref{sec:DivFreeHDG}.
As already explained there, an exact numerical integration of the polynomial integrand in this term requires a quadrature of order $3k+3$, which means order 27 for $k=8$.
If instead the quadrature is chosen to be only of order $2k=16$, already this seemingly small under-integration error accumulates over time to such an extent that the last pairing is impacted.
Figure \ref{fig:Palinstrophy_INACCURATE_RT8} (right) shows what happens with the computational results if exclusively in the convection term the quadrature of order $27$ is replaced by a quadrature of order $16$.
At least on the coarse meshes, the palinstrophy changes after $\tbar=260$, which again underlines the sensitivity of the problem.
On sufficiently fine meshes, this source of error fortunately does not seem to have a large impact anymore.

%--------------------------------------Compiler settings-----------------------------------------
\subsection{Different compiler settings (FMA)}

Lastly, let us document that even such a thing as the particular compiler setting can have an impact on the results computed for the Kelvin--Helmholtz instability problem.
In particular, we encountered inconsistencies in the results related to the question of whether fused multiply-add (FMA) is used or not.
Whenever an operation $a \leftarrow a+\rb{b\cdot c}$ is made two roundings are performed (noFMA).
When only one rounding is performed, on the other hand, the compiler setting FMA is used.
All results  presented in Section~\ref{sec:ComputationalStudies} were computed without FMA.

\begin{figure}[b]
	\centering
	\includegraphics[width=0.49\textwidth]{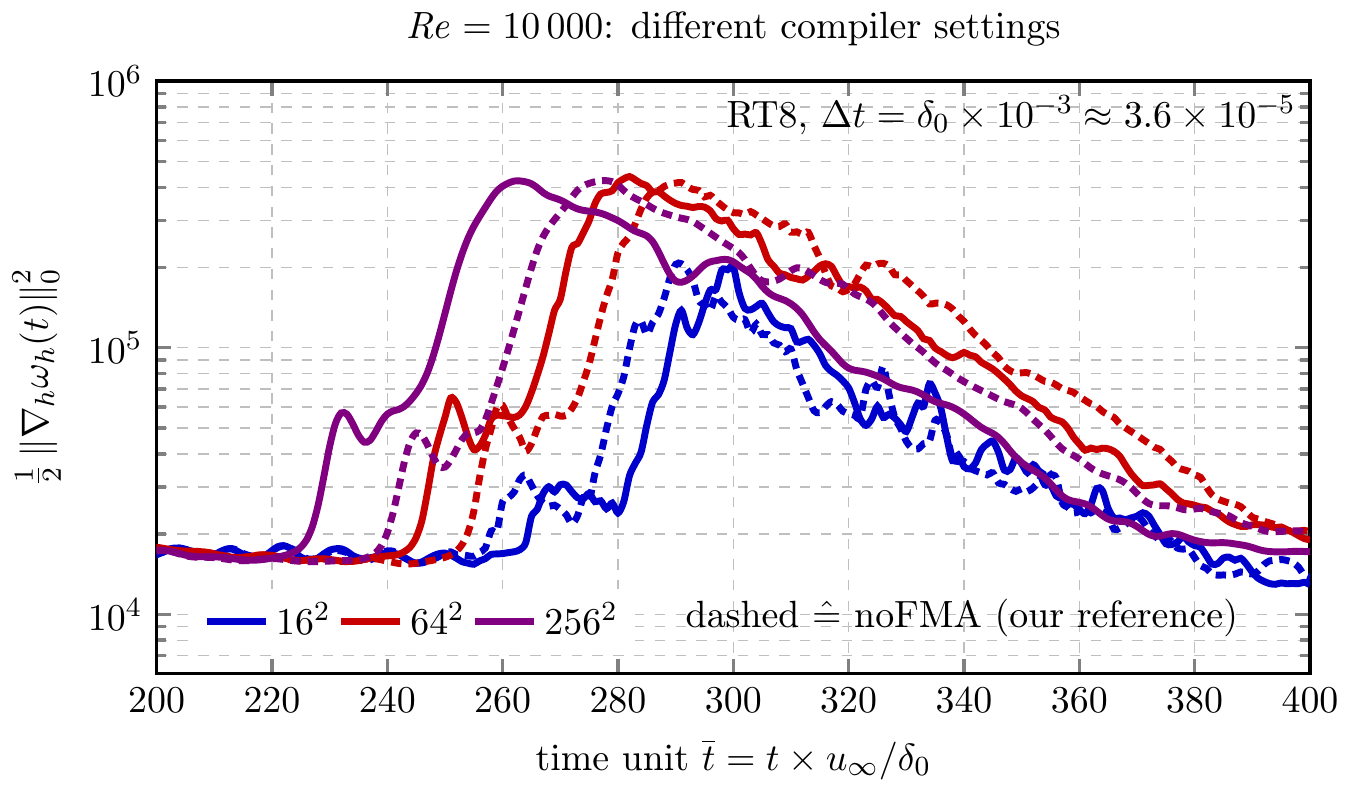} 
	\caption{Evolution of palinstrophy $\Pal\rb{t,\uu_h}$ for $\Rey=\num{10000}$ with RT8 HDG and SBDF2 IMEX time stepping on a sequence of square meshes. Perturbation due to FMA vs.\ noFMA compiler setting. Right: Dashed line shows the corresponding noFMA solutions from Section \ref{sec:ComputationalStudies}. For $\tbar < 200$, there is no visible difference. }
	\label{fig:Palinstrophy_COMPILER_RT8}
\end{figure}

Figure~\ref{fig:Palinstrophy_COMPILER_RT8} shows how the computational results change if the compiler setting noFMA is replaced by FMA.
The differences are not as extreme as for less accurate solutions of linear systems, for example. 
But still, these results underline how sensitive simulations can react to perturbations which are due solely to accumulated round-off errors.

%-----------------
%-----------------

%------------------------------------------------------------------------------------------------
%----------------------------------------CONCLUSIONS---------------------------------------------
%------------------------------------------------------------------------------------------------
\section{Summary and conclusions}
\label{sec:Conclusions}

In this work, we considered a 2D Kelvin--Helmholtz instability problem at various Reynolds numbers.
Section~\ref{sec:KHInstabilityProblem} dealt with the setting of the problem and an overview of corresponding frequently considered quantities of interest was given.
The most important quantity, in our opinion, turned out to be the palinstrophy \textemdash{} half the squared $\LP{2}{}$-norm of the gradient of the vorticity.
It is the most sensitive measure for indicating vortex merging in the considered problems.

From the theoretical point of view, in Section~\ref{sec:SelfOrga2DTurbulence}, the theory of self-organization in two-dimensional incompressible flows at high Reynolds numbers, originally presented by Van Groesen, was applied to the Kelvin--Helmholtz instability problem.
Here, the continuous problem was considered and the main new aspect laid in free-slip and periodic boundary conditions.
Again, the palinstrophy, along with kinetic energy and enstrophy, played a crucial role.
The most important result is Lemma~\ref{lem:InvariantSets}, which basically says that the considered problem is extremely sensitive with respect to (small) perturbations.

Section~\ref{sec:DivFreeHDG} explained both the spatial and temporal discretization schemes that were used to obtain numerical results.
A high-order ($k=8$), exactly divergence-free $\HDIV$-conforming Hybrid Discontinuous Galerkin (HDG) method in space and a multistep implicit-explicit time stepping scheme based on BDF2 were applied.

With the intent of keeping perturbations stemming from the discretization as small as possible, Section~\ref{sec:ComputationalStudies} presented computational studies of the Kelvin--Helmholtz instability problem for Reynolds numbers $\Rey\in\set{\num{100},\num{1000},\num{10000}}$ in a comparatively benign and controllable situation.
We showed that while it is possible to obtain mesh-converged simulations with respect to kinetic energy for all Reynolds numbers, perfectly controlling the enstrophy is already hard for high Reynolds numbers.
Concerning the palinstrophy, we obtained reliable reference results up to $\tbar=200$, which corresponds to a situation with two rotating vortices.
As it is manifesting in the palinstrophy, the point in time of the last merging, on the other hand, could not be controlled reliably.

In order to emphasize how sensitive the Kelvin--Helmholtz instability problem is, Section~\ref{sec:Perturbations} presented results for some common sources of perturbation: structured and unstructured triangular meshes, solvers of linear systems, numerical integration, and compiler settings.
All of them were studied to underline that the last merging process is extremely prone to perturbations of any kind.
In light of these results, it is not surprising that obtaining conclusive results with respect to all quantities of interest, over the whole time of the simulation, remains an open problem.

The data of the results of our computational studies are openly available for the community as explained in the Appendix.

%------------------------------------------------------------------------------------------------
%------------------------------------------------------------------------------------------------
\section{Appendix}
\label{sec:Appendix}

The detailed data of the computational results obtained in this study are presented at \url{https://ngsolve.org/kh-benchmark} \cite{KHpage}. 
Data which allows for a validation or comparison with other results can be accessed from there. 
This includes the data corresponding to all plots in this work, i.e.\ the time series for kinetic energy, enstrophy, vorticity thickness and palinstrophy, the energy spectra and field data for the velocity and vorticity at selected times.

%------------------------------------------------------------------------------------------------
%------------------------------------------------------------------------------------------------

%------------------------------------------------------------------------------------------------
%------------------------------------------------------------------------------------------------
\section*{Acknowledgments}
%------------------------------------------------------------------------------------------------
%------------------------------------------------------------------------------------------------
Philip L.\ Lederer has been funded by the Austrian Science Fund (FWF) through the research program ``Taming complexity in partial differential systems'' (F65) - project ``Automated discretization in multiphysics'' (P10). \\
We gratefully acknowledge the comments and suggestions on this manuscript from the anonymous reviewers.

%------------------------------------------------------------------------------------------------
%------------------------------------------------------------------------------------------------
\def\bibsection{\section*{References}}
\addcontentsline{toc}{section}{References}
%------------------------------------------------------------------------------------------------
%------------------------------------------------------------------------------------------------
\bibliography{KH-benchmark-BibTeX.bib}
\bibliographystyle{abbrvurl}

\end{document}